\documentclass[12pt]{article}
\RequirePackage[OT1]{fontenc}
\RequirePackage{amsthm,amsmath}

\RequirePackage[colorlinks,citecolor=blue,urlcolor=blue]{hyperref}

\usepackage{graphicx}
\usepackage{amsmath}
\usepackage{graphics}
\usepackage{epsfig}
\usepackage{amssymb}
\usepackage{color}

\usepackage[authoryear]{natbib}

\newcommand{\T}{\mathcal{T}}

\newcommand{\be}{\begin{eqnarray}}
\newcommand{\ee}{\end{eqnarray}}
\newcommand{\bea}{\begin{eqnarray*}}
\newcommand{\eea}{\end{eqnarray*}}

\newcommand{\ve}{\varepsilon}

\newtheorem{theorem}{Theorem}[section]
\newtheorem{lemma}{Lemma}[section]
\newtheorem{proposition}{Proposition}[section]

\newtheorem{example}{Example}[section]

\newtheorem{remark}{Remark}[section]


\begin{document}

\title{Prediction in regression models with continuous observations}

 \author{Holger Dette$^{1,4}$, Andrey Pepelyshev$^2$, Anatoly Zhigljavsky$^3$}

\maketitle

 \footnotetext[1]{Fakult\"at f\"ur Mathematik, Ruhr-Universit\"at Bochum,
Bochum, 44780, Germany. holger.dette@rub.de}
 \footnotetext[2]{School of Mathematics, Cardiff  University, Cardiff, CF24 4AG, UK. pepelyshevan@cardiff.ac.uk}
 \footnotetext[3]{School of Mathematics, Cardiff  University, Cardiff, CF24 4AG, UK. ZhigljavskyAA@cf.ac.uk}
 \footnotetext[4]{Corresponding author.}

\begin{abstract}

We consider the problem of predicting  values of a random process or field  satisfying  a linear  model $y(x)=\theta^\top f(x) + \ve(x)$, where  errors $\ve(x)$ are correlated. This is a common problem  in
 kriging, where the case of discrete observations is standard. By  focussing on the case of continuous observations, we  derive expressions for the best linear unbiased predictors and their mean squared error.
 Our  results are also applicable in the case where the derivatives of the process $y$ are  available,
and either a response or one of its  derivatives need to be  predicted.  The theoretical results are illustrated
by several  examples in particular for the popular Mat\'{e}rn  $3/2$ kernel.

\end{abstract}





\section{Introduction}
\def\theequation{1.\arabic{equation}}
\setcounter{equation}{0}

A common problem, which  occurs in many different areas, most notably geostatistics \citep{ripley1991statistical,cressie1993statistics}, computer experiments
\citep{sacks1989design, stein1999interpolation, santner2003design,leatherman2017designing} and machine learning \citep{rasmussen2006gaussian}, is
to predict the response $y(t_0)$ at a point $t_0\in\mathbb{R}^d$ from
given responses $y(t_1),\ldots,y(t_N)$ at points $t_1,\ldots,t_N\in\mathbb{R}^d$, where $t_0\neq t_i$ for all $i=1,\ldots,N$.
Making the prediction assuming that responses are observations of a random field is called kriging  \citep{stein1999interpolation}.
In classical kriging, it is assumed that $y$ is a random field of the form
\be
\label{1.1}
 y(t)=f^\top(t)\theta+\epsilon(t),
\ee
where $f(t) \in \mathbb{R}^m $ is a vector of known regression functions, $\theta \in \mathbb{R}^m $ is a vector of unknown parameters and
$\epsilon$ is a random field with zero mean and  existing covariance kernel, say $K(t,s) = E [\epsilon(t)\epsilon(s) ]$. The components of the vector-function $f(t)$ are assumed to be linearly independent on the set of points where the observations have been made.

It is well-known, see e.g. \cite{sacks1989design}, that in the case of discrete observations the best linear unbiased predictor (BLUP) of $y(t_0)$ has the form
\be
 \label{eq:kriging}
 \hat y(t_0)=f^\top(t_0)\hat\theta_{\mathrm{BLUE}}+K^\top_{t_0}\Sigma^{-1}   (Y-X\hat\theta_{\mathrm{BLUE}} ),
\ee
where  $\Sigma=\big(K(t_i,t_j)\big)_{i,j=1}^N$
is an $N\!\times\! N$-matrix, $K_{t_0}=\big(K(t_0,t_1),\ldots,K(t_0,t_N)\big)^\top$ is a vector in $\mathbb{R}^N$,
$X=(f(t_1),\ldots,f(t_N))^\top $ is an $N\!\times\! m$-matrix, $Y=(y(t_1),\ldots,$ $y(t_N))^\top  \in \mathbb{R}^N$ is a vector of observations and
$$\hat\theta_{\mathrm{BLUE}}=(X^\top \Sigma^{-1}X)^{-1}X^\top \Sigma^{-1}Y$$ is the  best linear unbiased estimator (BLUE) of $\theta$.
The BLUP satisfies the unbiased condition
$
 \mathbb{E}[\hat y(t_0)]=\mathbb{E}[y(t_0)]
$
and minimizes the mean squared error  $\mathrm{MSE}(\tilde y(t_0))=\mathbb{E}\left(y(t_0) -\tilde y(t_0)\right)^2$
in the class of all linear unbiased predictors $\tilde y(t_0)$; its mean squared error is
\bea
 \mathrm{MSE}(\hat y(t_0))=K(t_0,t_0)-\left[\begin{matrix}f(t_0)\\ K_{t_0}\end{matrix}\right]^\top
 \left[\begin{matrix}0&X^\top \\X&\Sigma \end{matrix}\right]^{-1}
 \left[\begin{matrix}f(t_0)\\ K_{t_0}\end{matrix}\right].
\eea

In the present paper, we generalize the predictor \eqref{eq:kriging} to the case of
continuous observations of the response including possibly derivatives and { prediction of derivatives and weighted averages of $y(t)$. We shall separately  consider the cases where the observation region is an interval or a product set (in particular, square).  }

An important observation concerning construction of the BLUPs at different points is the fact that there is a considerable common part related to the use of the same BLUE. This could lead to significant computational savings relative to independent construction of the BLUPs. This observation extend to the cases when the observations are taken in $\mathbb{R}^d$ and when derivatives are also used for predictions.

The remaining part of this paper is organized as follows. In Section  \ref{sec2} we consider the BLUPs when we observe the process or field only.
In Section \ref{sec3} we study the BLUPs for either process values or one of its derivatives when we observe the process (or field)  with derivatives.
In Section \ref{sec4}  we provide  proofs of the main results and in  an Appendix we give more illustrating examples of the BLUPs for  particular kernels.

\section{ Prediction without derivatives } \label{sec2}
\def\theequation{2.\arabic{equation}}
\setcounter{equation}{0}

\subsection{Prediction  at a point}

Assume $\T\subset \mathbb{R}^d$ and
consider  prediction  at a point $t_0\not\in\T$ for a response given by the model \eqref{1.1}, where the  observations for all $t\in\T$ are available.
 The vector-function  \mbox{$f\!:  \T \!\to \! \mathbb{R}^m $} is assumed to contain functions which are bounded,  integrable, smooth enough and linearly independent on $\T$;  the covariance kernel $K(t,s) = E[\epsilon (t) \epsilon (s)  ]$  is assumed strictly positive definite.

A general linear  predictor of $y(t_0)$ can be defined as
$$
 \hat y_Q(t_0)=  \int_\T y(t)Q(dt),
$$
where  $Q$ is a signed measure
defined on the Borel field of~$\T$.
This predictor is unbiased if $\mathbb{E}[\hat y_Q(t_0)]=\mathbb{E}[y(t_0)]$, which is equivalent to the condition
$$
\int_\T f(t)Q(dt)= f(t_0).
$$

The mean squared error (MSE) of $\hat y_Q(t_0)$ is given by
$$
\mathrm{MSE}(\hat y_{Q}(t_0))=\mathbb{E}\left(y(t_0) -\hat y_{Q}(t_0)\right)^2 \, .
$$
The best linear unbiased predictor (BLUP) $\hat y_{Q_*}(t_0)$ of $y(t_0)$  minimizes  the mean squared error
$\mathrm{MSE}(\hat y_{Q}(t_0))$ in the set of all linear unbiased predictors.
The corresponding signed measure $Q_*$  will be called BLUP measure throughout this paper.
Unlike  the case of discrete observations, the BLUP measure does not have to exist for continuous observations.

\medskip

{\bf Assumption A.} {\it

(1) The best linear unbiased estimator (BLUE)
  $\hat{\theta}_{\mathrm{BLUE}} = \int_\T y(t) G(dt) $  exists in the model \eqref{1.1}, where $G(dt)$ is some  signed vector-measure on $\T$, \\
 (2) There exists
 a signed measure $\zeta_{t_0}(dt)$ which satisfies  the equation
\be
\label{eq:zeta}
 \int_\T  K(t,s)\zeta_{t_0}(dt)= K(t_0,s),\;\; \forall s\in \T.
\ee
}

Assumption A will be discussed in Section \ref{sec22}  below. We  continue  with
a general statement establishing the existence and explicit form of the BLUP.

\begin{theorem}
\label{th:predict-t0}
If  Assumption A holds
then
 the  BLUP measure $Q_*$ exists and  is given by
\be
\label{eq:blupA}
Q_*(dt)=\zeta_{t_0}(dt)+c^\top G(dt),
\ee
where the  signed measure $\zeta_{t_0}(dt)$
 satisfies \eqref{eq:zeta}
and
$
c= f(t_0)-\int_\T f(t)\zeta_{t_0}(dt)\, .
$
The MSE  of the corresponding BLUP  $\hat{y}_{Q_*}(t_0)$ is given by
\be
\label{eq:blup_mse0}
\mathrm{MSE}(\hat{y}_{Q_*}(t_0))=K(t_0,t_0)+c^\top  D f(t_0)- \int_{\T} K(t,t_0) Q_*(dt)\, ,
\ee
where $
D= \int_\T\!\! \int_\T K(t,s) G(dt) G^\top (ds)
$
 is the covariance matrix of  $\hat{\theta}_{\mathrm{BLUE}}=\int_\T y(t) G(dt) $.
\end{theorem}

This theorem is a particular case of a more general Theorem \ref{th:predict-nu}, which considers the problem of predicting an integral of the response.
A few examples illustrating applications of Theorem~\ref{th:predict-t0} for particular kernels are given in the  Appendix.

We can interpret the construction  of the BLUP at $t_0$ in  model  \eqref{1.1} as the following two-stage algorithm.
At stage 1, we use  the BLUE
$\hat{\theta}_{\mathrm{BLUE}}=\int_\T y(t) G(dt) $
for estimating  $\theta$. At stage 2, we compute the BLUP  in the model
$$
\tilde{y}(t)=y(t)- f^\top (t) \widehat{\theta}_{BLUE} = {\varepsilon}(t)-f^\top (t){\int_\T {\varepsilon}(t') G(dt')}\, ,
$$
which is a model with new error process and no trend. Straightforwardly,  the covariance function of the process $\tilde{y}(t)$ is
calculated as
$$
\tilde{K}(t,s)= K(t,s)-f^\top(t) D f(s).
$$
It then follows from Theorem~\ref{th:predict-t0} applied to the new model  that the signed measure $Q_*(dt)$ satisfies the equation
\bea
\label{eq:zeta1}
 \int_\T  \tilde{K}(t,s)Q_*(dt)= \tilde{K}(s,t_0),\;\; \forall s\in \T\, .
\eea
From \eqref{eq:blup_mse0} in the new model, we obtain an alternative representation  for the MSE of the BLUP  $\hat{y}_{Q_*}(t_0)$; that is,
\bea
\mathrm{MSE}(\hat{y}_{Q_*}(t_0))=\tilde{K}(t_0,t_0)- \int_{\T} \tilde{K}(t,t_0) Q_*(dt)\, .
\eea

\subsection{Validity of Assumption A} \label{sec22}

If  ${\cal T}$ is a discrete set then  Assumption A  is satisfied for any strictly positive definite covariance kernel.

In general, the main part of Assumption A is the existence of the BLUE of the parameter  $\theta$, which has been clarified by
 \cite{DPZ2018}. According to their Theorem 2.2, the BLUE of $\theta$ exists if and only if
there exists   a signed vector-measure $G=(G_1, \ldots, G_m)^\top $ on $\T$, such that the $m\!\times\! m $-matrix $\int_\T f(t) G^\top (dt)$ is the identity matrix and
\be
  \int_\T  K(t,s)G(dt)= D f(s)
 \label{eq:suff-ness-cond-qG}
\ee
holds for  all $s\in \T$ and  some $m\! \times\! m$-matrix $D$. In this case,
$
\widehat{\theta}_{BLUE} = \int_\T Y(t)G(dt)
$
and
$
D
$
is
the covariance matrix of
 $\widehat{\theta}_{BLUE}
  $; this matrix  does not have to be non-degenerate.

Let $\cal{H}_K$ be the reproducing kernel Hilbert space (RKHS) associated with kernel $K$.
If the function $K(t_0,s)$ belongs to $\cal{H}_K$, then the second part of Assumption A is also satisfied; that is, there exists a
measure $\zeta_{t_0}(dt)$ satisfying the equation \eqref{eq:zeta}. This follows from results of \cite{parzen1961approach}. Note that
the function $K(t_0,s)$ does not automatically belong to $\cal{H}_K$ since in general  $t_0 \notin  \T$.

If
all components of $f$ belong to $\cal{H}_K$   then Assumption A holds and
the matrix $D$ in \eqref{eq:suff-ness-cond-qG} is non-degenerate; see \cite{DPZ2018} and \cite{parzen1961approach}.

If the matrix $D$ in Theorem \ref{th:predict-t0} is non-degenerate then this theorem can be reformulated in the following form which is practically more convenient
as there is no unbiasedness condition to check.

\begin{proposition}
\label{prop1}
Assume that there exists
 a signed measure $\zeta_{t_0}(dt)$ satisfying \eqref{eq:zeta}
and
a signed vector-measure
$\zeta(dt)$  satisfying equation
\be
\label{eq:zeta0}
 \int_\T  K(t,s)\zeta(dt)= f(s),\;\; \forall s\in \T.
\ee
If additionally  the matrix
$
C=\int_\T f(t)\zeta^\top (dt)
$
is non-degenerate, then the  BLUP measure exists and  is given by
 \eqref{eq:blupA} with $D=C^{-1}$. Its MSE is given by  \eqref{eq:blup_mse0}.
\end{proposition}

Clearly, if the conditions of Proposition~\ref{prop1} are satisfied then the BLUE measure $G(dt)$ is expressed via the measure $\zeta(dt)$ by $G(dt)=C^{-1}\zeta(dt)$.


Explicit forms of the BLUP for some kernels are given in the  Appendix.




\subsection{Matching expressions in the case of discrete observations}

Let us show that in the case of discrete observations the form of the BLUP of Proposition \ref{prop1} coincides with the standard form \eqref{eq:kriging}.
Assume  that ${\cal T}$ is finite, say, ${\cal T}=\{t_1,\ldots,t_N\}$.
In this case, equation \eqref{eq:zeta0} has the form $ \Sigma \zeta = X$, where $\zeta$ is and
$N\!\times\! m$-matrix. {Since the kernel $K$ is strictly positive definite, this} gives
  $\zeta=\Sigma^{-1}X$, and we also obtain
  $C=X^\top \Sigma^{-1}X$, $G^\top =C^{-1}\zeta^\top $. A general linear predictor is of form $\tilde{y}(t_0)=Q^\top Y$ and the BLUP is $Q_*^\top Y$ with
$
 \label{eq:discrBLUP}
Q^\top _*=\zeta_{t_0}^\top +c^\top G^\top ,
$
where $\zeta_{t_0}=\Sigma^{-1}K_{t_0}$
satisfies equation \eqref{eq:zeta}
and
$
 c= f(t_0)-X^\top \zeta_{t_0}.
$
Expanding the expression for $Q^\top _*$ we obtain
\be
\label{eq:q-discr}
 Q^\top _*&=&(\Sigma^{-1}K_{t_0})^\top +c^\top C^{-1}(\Sigma^{-1}X)^\top  \nonumber\\
 &=&K_{t_0}^\top \Sigma^{-1}+(f(t_0)-X^\top \Sigma^{-1}K_{t_0})^\top C^{-1}X^\top \Sigma^{-1}\, .
\ee

The classical form  of the BLUP is given by  \eqref{eq:kriging}, which can be written as $Q^\top  Y$ with
$
 Q^\top =f^\top (t_0)C^{-1}X^\top \Sigma^{-1}+K^\top _{t_0}\Sigma^{-1}-K^\top _{t_0}\Sigma^{-1}XC^{-1}\Sigma^{-1}X^\top .
$
and coincides with \eqref{eq:q-discr}.

\subsection{Predicting an average with respect to a measure}

Assume that we have a realization of a random field  \eqref{1.1}
observed for all $t\in\T \subset \mathbb{R}^d$.
Consider the prediction problem of
$Z=\int_{\cal S} y(t) \nu(dt),$
where $\nu(dt) $ is some (signed) measure on the Borel  field  of $\mathbb{R}^d$ with support~${\cal S}$. Assume that ${\cal S} \setminus \T \neq \emptyset$ (otherwise, if ${\cal S} \subseteq \T$, the problem is trivial as
we observe the full trajectory $\{y(t) ~|~ t \in \T \} $).
We interpret $Z$ as a weighted average of the true process values on~${\cal S}$.
The  general linear  predictor can be defined as
\begin{align}\label{zqaht}
 \hat Z_Q=  \int_\T y(t)Q(dt),
\end{align}
where  $Q$ is a signed measure
on the Borel  field  of $\T$. The estimator $\hat Z_Q$ is unbiased if
and only if
\be
\label{2.3.1}
\int_\T f(t)Q(dt)= \int_{\cal S} f(s) \nu(ds)\, .
\ee

The  BLUP signed measure  $Q_*$ minimizes
$$
\mathrm{MSE}(\hat Z_{Q})=\mathbb{E}\big(Z -\hat Z_{Q}\big)^2
$$
among all signed measure $Q$ satisfying the unbiasedness condition \eqref{2.3.1}.
Assumption A and Theorem \ref{th:predict-t0} generalize to the following.

\medskip

{\bf Assumption A$^\prime$.} {\it
The BLUE
  $\hat{\theta}_{\mathrm{BLUE}}$  exists
 and there exists
 a signed measure $\zeta_{\nu}(dt)$ which satisfies  the equation
\be
\label{eq:zeta1}
 \int_\T  K(t,s)\zeta_{\nu}(dt)= \int_{{\cal S}} K(s,u)\nu(du),\;\; \forall s\in \T.
\ee
}

\begin{theorem}
\label{th:predict-nu}
Suppose that
 Assumption A$^\prime$ holds and
let $D$ be the covariance matrix of $\hat{\theta}_{\mathrm{BLUE}}=\int_\T y(t) G(dt) $.
Then the  BLUP measure exists and  is given by
\be
\label{eq:blup1}
Q_*(dt)=\zeta_{\nu}(dt)+c^\top G(dt),
\ee
where $\zeta_{\nu}(dt)$
is the  signed measure satisfying \eqref{eq:zeta1}
and
$
c= 
\int_{{\cal S}}f(s)\nu(ds)-\int_\T f(t)\zeta_{\nu}(dt) 
\, .
$
The MSE  of the BLUP  $\hat Z_{Q_*}$ is given by
\bea
\mathrm{MSE}(\hat Z_{Q_*})\!=\!\!\int_{{\cal S}} \!\int_{{\cal S}}\! K(s,u)\nu(\!ds\!)\nu(\!du\!)\!+\!c^\top \! D\!\!\int_{\cal S}\!f(s) \nu(\!ds\!)\!-\! \int_{\cal S}\!\int_{\T}\! K(t,u) \nu(\!du\!) Q_*(\!dt\!) .
\eea
\end{theorem}

The proof of Theorem \ref{th:predict-nu} is given in Section \ref{sec4} and contains the proof of Theorem \ref{th:predict-t0} as a special case.
Note also that the BLUP $\hat Z_{Q_*}$ is simply the average (with respect to the measure $\nu$) of the BLUPs at points $s \in {\cal S}$.

\subsection{Location scale model on a product set} \label{sec24}


In this section we consider the location scale  model
\be \label{eq:loc_model}
y(t)=\theta+\ve (t),
\;\;{\rm  where}\; t=(t_1,t_2)\in  \T
\ee
 and assume that
the kernel $ \mbox{\textsf{K}}$ of the random field ${\ve}(t)$ is given by
\be
\label{eq:loc_prod}
\mbox{\textsf{K}}( t,t')&=&\mathbb{E}[{\ve}(t){\ve}(t')]=K_1(t_1,t_1')K_2(t_2,t_2')~,~
\ee
for $t=(t_1,t_2), t'=(t'_1,t'_2)\in \T.$ We also assume that the set $\T \subset \mathbb{R}^2$ is a product-set of the form  $\T = \T_1 \times  \T_2$,
where $\T_1 $ and $ \T_2$ are Borel subsets of~$\mathbb{R}$ (in particular, these sets could be discrete or continuous).
{The kernel $K$ of the product form  \eqref{eq:loc_prod} is  called separable; such kernels are frequently used in    modelling of
 spatial-temporal structures because they  offer enormous computational benefits, including rapid  fitting
 and simple extensions of many techniques  from time series and classical geostatistics [see \cite{Gneitingetal2006} or   \cite{fuentes2006} among many others].}

Assume that Assumption A$^\prime$ holds for two one-dimensional models
\be
\label{eq:prod_models}
y_{(i)}(u)=\theta + \ve_{(i)}(u)~,~~u \in \T_i ~~~~ (i=1,2)
\ee
with
$
\nonumber
K_i(u,u')=\mathbb{E}[\ve_{(i)}(u)\ve_{(i)}(u')] ~,~~u,u' \in \T_i ~~ (i=1,2).
$
Let    the measures $G_i(du)$   define  the BLUE
$
\int_{\T_i} y_{(i)}(u) G_i(du)
$
in these two models.
Then the BLUE of $\theta$ in the model \eqref{eq:loc_model}  is given by
$
\hat {\theta}=
\int_\T y(t)  \mbox{\textsf{G}}(dt),
$
where $  \mbox{\textsf{G}}$ is a product-measure
$
 \mbox{\textsf{G}}(dt)=G_1(dt_1)G_2(dt_2),
$
Assume we want to predict $ y(t)$ at a point $T=(T_1,T_2) \notin \T$.
Note that  equation \eqref{eq:zeta1}  can be rewritten as
\bea
\label{eq:zetaZZ}
 \int_\T  \mbox{\textsf{K}}(t,s) \zeta_{T} (dt)
  =   \mbox{\textsf{K}}(s,T),\;\; \forall s\in \T.
\eea

A solution of the above equation has the form
$
 {\zeta}_{T}(dt_1,dt_2)= \zeta_{T_1}(dt_1) \zeta_{T_2}(dt_2)\, ,
$
where $\zeta_{T_i}(dt)$ ($i=1,2$) satisfies the equation
\be
\label{eq:zetaZZ1}
 \int_{\T_i}  K_i(u,v) \zeta_{T_i}(du)=  K_i(v,T_i) ,\;\; \forall v\in \T_i.
\ee

Finally, the BLUP at the point $T=(T_1,T_2)$ is $\int_\T y(t)  \textsf{Q}_{*}(dt) $,
where
\bea
 \textsf{Q}_{*}(dt)=  \zeta_{T} (dt)+ c \, \mbox{\textsf{G}}(dt)
\;\; {\rm
with } \;\;
c=1-  \int_\T  \zeta_{T} (dt).
\eea

The measure $ \mbox{\textsf{G}}(dt)$ is the BLUE measure and does not depend on $T_1,T_2$. On the other hand, the measure $ {\zeta}_{T}(dt)$ and constant $c$ do depend on $T_1,T_2$.
The MSE of the BLUP is $\mathrm{MSE}(\hat{ y}_{Q_*}(T))=1+c-\int_\T  \mbox{\textsf{K}} (t,T) \textsf{Q}_{*}(dt)\, .$

\begin{example}
{\rm
Consider the case of $\T=[0,1]^2$ and the exponential kernel
\bea
 \mbox{\textsf{K}}(t,t')&=&\mathbb{E}[ \ve(t) \ve(t')]=\exp \left\{-\lambda \left[|t_1-t_1'|+|t_2-t_2'| \right] \right \},\;\;
\eea
where $\lambda>0$ and $t=(t_1,t_2), t'=(t'_1,t'_2)
 \in [0,1]^{2}$.
Define the measure
$$
G(du) = \frac1{2+\lambda} \left[  \delta_0 (du)+ \delta_1 (du)+\lambda  du \right], \;\; u \in [0,1].
$$
In view of \citep[Sect 3.4]{DPZ2018},  $\int_0^1 y(u) G(du)$ is the BLUE in the model $y(u)=\theta + \ve(u)$ with kernel $K(u,u')=\mathbb{E}[\ve(u)\ve(u')]=e^{-\lambda|u-u'|}$, $u,u' \in [0,1]$.
The equation \eqref{eq:zetaZZ1} can be rewritten as
\bea
\label{eq:zetaZZ2}
 \int_0^1  e^{-\lambda|v-u|}  \zeta_{T_i}(du)=  e^{-\lambda|v-T_i|} ,\;\; \forall v\in [0,1].
\eea
It follows from \citep[Sect 3.4]{DPZ2018} that this equation is satisfied by the measure
\bea
 \zeta_{T_i}(du)= \left\{
                    \begin{array}{ll}
                     e^{-\lambda |T_i|} \delta_0(du),  &  \mbox{ if ~} T_i \leq 0, \\
                      \delta_{T_i}(du),   &   \mbox{ if ~}  0 \leq T_i \leq 1,  \\
                      e^{-\lambda (T_i-1)} \delta_1(du), &   \mbox{ if ~}  T_i \geq 1. \\
                    \end{array}
                  \right.
\eea

For   $T_1 \leq 0$ we obtain $ \textsf{Q}_{*}(dt)= {\zeta}_{(T_1,T_2)}(dt)+ c  \mbox{\textsf{G}}(dt)$ in the following form
\bea
  \textsf{Q}_{*}(dt)  = \left\{
                    \begin{array}{l}
                     e^{-\lambda |T_1|}\delta_0(dt_1)\delta_{T_2}(dt_2) +\left(1- e^{-\lambda |T_1|} \right)   \mbox{\textsf{G}}(dt),
                         \mbox{~ if ~}  0 \leq T_2 \leq 1,  \\
                  e^{-\lambda |T_1|-\lambda |T_2|} \delta_0(dt_1) \delta_0(dt_2)       +\left(1- e^{-\lambda |T_1|-\lambda |T_2|} \right) \mbox{\textsf{G}}(dt),  \\
                      \mbox{~~~~~~~~~~~~~~~~~~~~~~~~~~~~~~~~~~~~~~~~~~~~~~~~~~~~~~~~~~~ if ~} T_2 \leq 0, \\
                      e^{-\lambda |T_1|-\lambda (T_2-1)} \delta_0(dt_1)\delta_1(dt_2) +\left(1- e^{-\lambda |T_1|-\lambda |T_2-1|} \right)  \mbox{\textsf{G}}(dt), \\
                         \mbox{~~~~~~~~~~~~~~~~~~~~~~~~~~~~~~~~~~~~~~~~~~~~~~~~~~~~~~~~~~~ if ~} T_2 \geq 1.
                    \end{array}
                  \right.
\eea
Similar formulas can be obtained for $0<T_1<1$ and $T_1\geq 1$. \\

In Table \ref{tab:rmse-exp-square} we show the square root of the MSE of the BLUP
for the equidistant design supported at points $(i/(N-1),j/(N-1))$, ${i,j=0,1,\ldots,N-1}$.
We can see that the MSE for the design with $N=4$ is already rather  close to the MSE for the design with large $N$ and the design with continuous observations.

\begin{table}[!hhh]
\caption{\it The square root of the MSE of the BLUP at several points for
the $N\!\times\! N$-point equidistant design in the location scale model on the square $[0,1]^2$
and the exponential kernel with $\lambda=2$. In the case $N=\infty$ we provide the  MSE for continuous observations.}
\label{tab:rmse-exp-square}
\begin{center}
\begin{tabular}{|c|c|c|c|c|c|c|c|}
\hline
 $N$ & 2 & 3 & 4 & 8 & 16  & 32 & $\infty$\\
\hline
$T=(2,2)$  &1.1446& 1.1225 & 1.1177 &1.1145& 1.11398 &  1.11386& 1.11383\\
$T=(0.5,2)$&1.1242& 1.0879 & 1.0884 &1.0831& 1.08177 &  1.08133& 1.08117\\
\hline
\end{tabular}
\end{center}
\end{table}

In Figure \ref{fig:rmse-ou} we show the plot of the square root of the MSE as a function of a prediction point for points $(T_1,T_2) \in [0.5,2]\times [0.5,2]$. As the design is symmetric with respect
to the   point $(0.5,0.5)$, the plot of the MSE is also  symmetric with respect to  this point. Consequently only the upper quadrant is depicted in the figure.

We observe that  the MSE tends to zero when the prediction point tends to one of design points and the MSE is almost constant  if the prediction point is far enough from  the observation domain.

\begin{figure}[!hhh]
\begin{center}
 \includegraphics[width=0.48\textwidth]{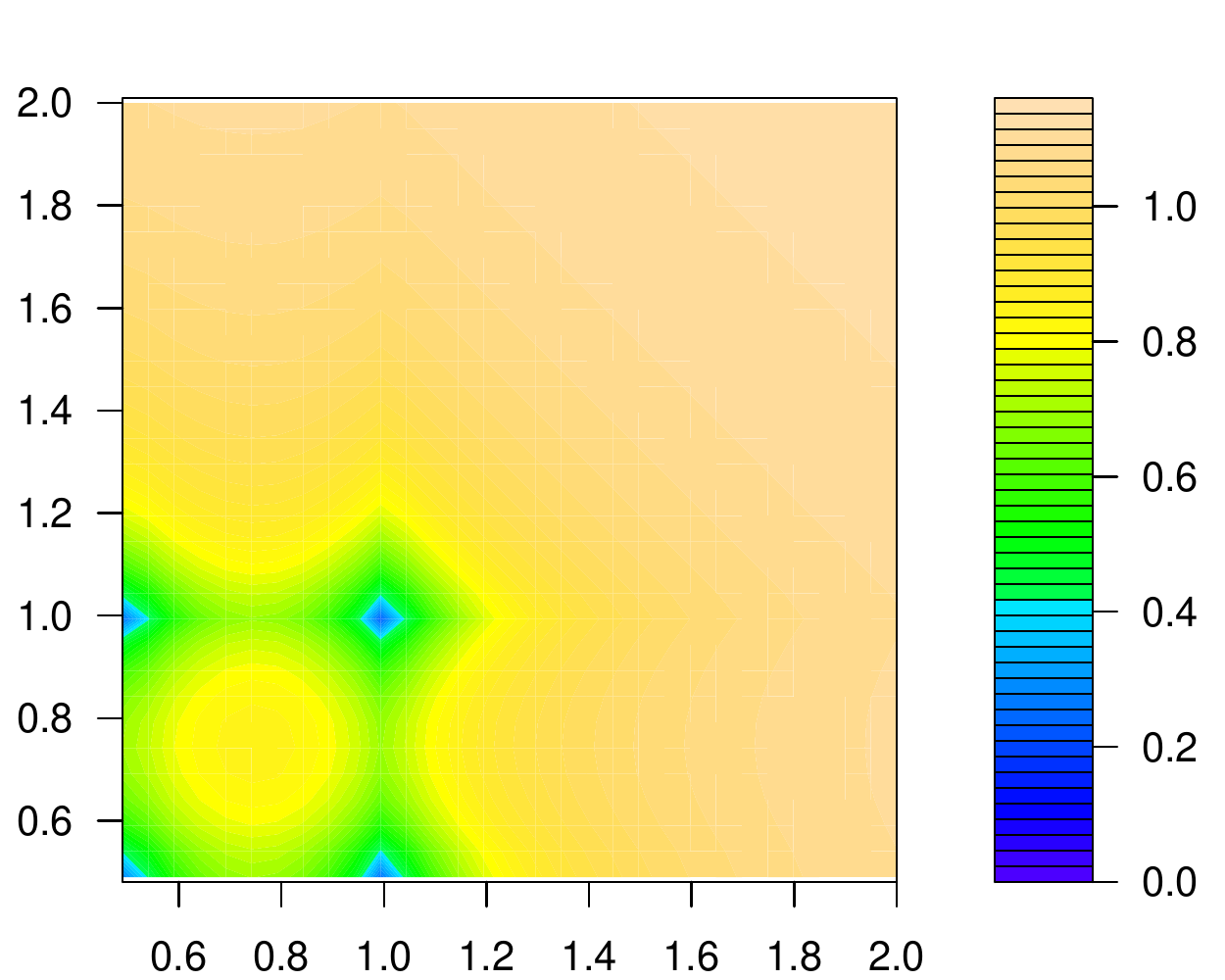}
 \includegraphics[width=0.48\textwidth]{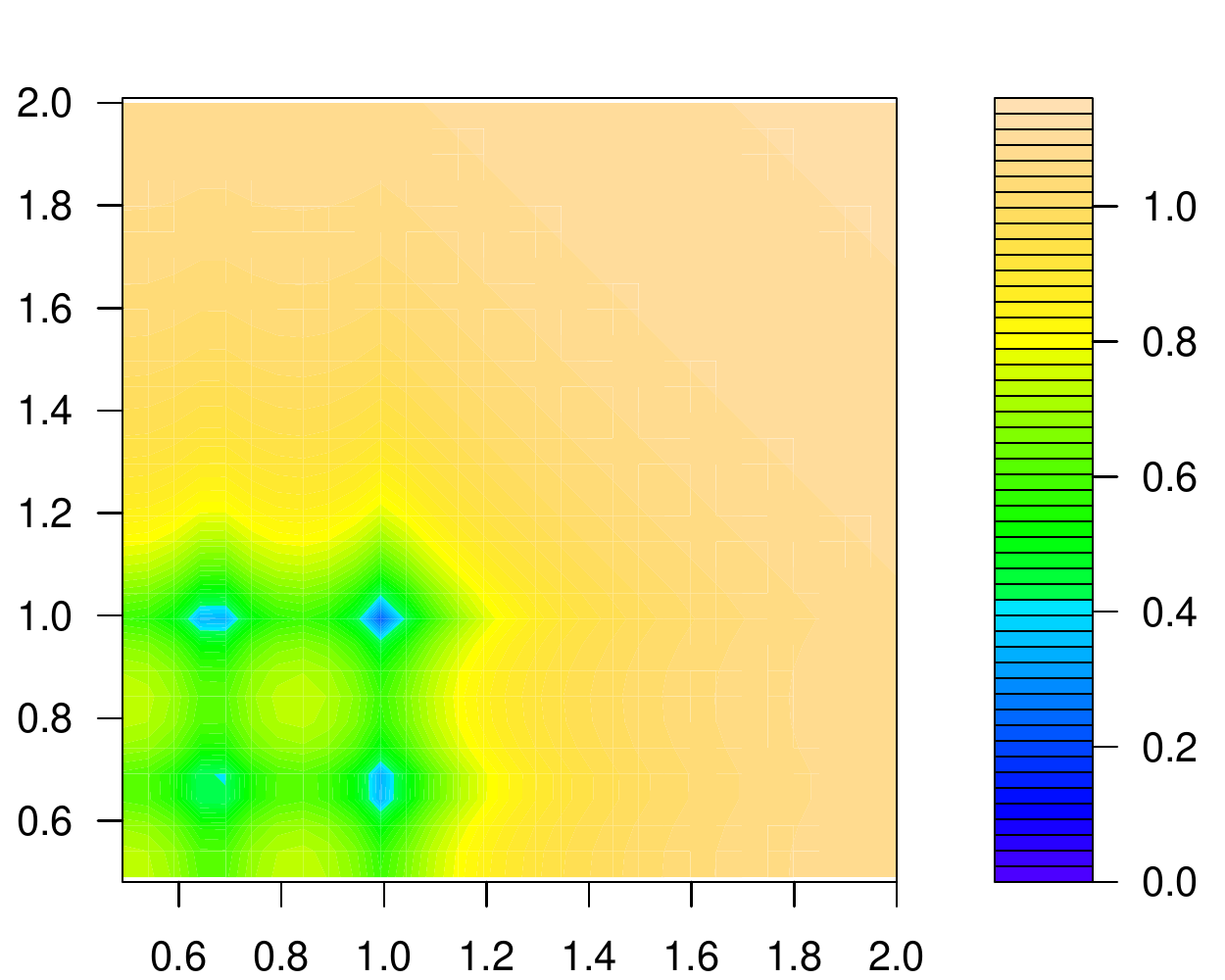}
\end{center}
\caption{\it The square root of the MSE of the BLUP for the $N\!\times\! N$-point equidistant design with $N=3$ (left) and $N=4$ (right) and the exponential kernel with $\lambda=2$.}
\label{fig:rmse-ou}
\end{figure}

}
\end{example}

\begin{remark} \label{rem1}
{\rm
The results of this section can be easily generalized to the case of $d>2$ variables and, moreover, to the model $y(t)=\theta f(t)+\ve (t)$, where $t=(t_1,\ldots,t_d)\in \T_1 \times \ldots \times \T_d$,
$
\mbox{\textsf{K}}( t,t')=\mathbb{E}[{\ve}(t){\ve}(t')]=K_1(t_1,t_1')\cdots K_d(t_d,t_d')
$
and
$
f(t)=f_{(1)}(t_1) \cdots f_{(d)}(t_d) ,
$ where $f_{(i)}$ are some functions on $\T_i~;~~i=1, \ldots, d$.
}
\end{remark}

\section{Prediction with derivatives} \label{sec3}
\def\theequation{3.\arabic{equation}}
\setcounter{equation}{0}

In this section we consider prediction problems, where the trajectory $y$ in model \eqref{1.1}  is differentiable (in the mean-square sense) and derivatives of the process (or field) $y$ are available. In Section~\ref{sec31} we discuss the discrete case of a once-differentiable process and in Section~\ref{sec32} we consider the general case of a $q$ times differentiable (in the mean-square sense) process  $y$ satisfying the model \eqref{1.1}. For the process $y$  to be $q$ times differentiable, the covariance kernel $K$ and vector-function $f$ in \eqref{1.1} have  to
be $q$ times differentiable, which is one of the assumptions in Section~\ref{sec32}.
In Section \ref{sec:sec33} we consider the prediction problem for the location scale model on a two-dimensional product set in the case where
 the kernel $ \mbox{\textsf{K}}$ of the random field ${\ve}$  has the product form
\eqref{eq:loc_prod}.
The results of this section  can be easily generalized to the case of $d>2$ variables.

\subsection{Discrete case} \label{sec31}

Consider the model \eqref{1.1}, where the kernel $K$ and vector-function $f$ are differentiable
 and one can observe  the process $y$ and its derivative at $N$ different points $t_{1} , \ldots , t_{N} \in \mathbb{R}$. In this case,
the BLUP of $y(t_0)$ has the form
\be
 \label{eq:kriging-deriv}
 \hat y(t_0)=f^\top (t_0)\hat\theta_{\mathrm{BLUE,}{2N}}+K^\top _{t_0,{2N}}\Sigma^{-1}   (Y_{2N}-X_{2N}\hat\theta_{\mathrm{BLUE,}{2N}} ),
\ee
where  $Y_{2N}=(y(t_1),\ldots,y(t_N),y'(t_1),\ldots,y'(t_N))^\top  \in \mathbb{R}^{2N}$,
$$
\Sigma=\left( {\Sigma _{00}\atop \Sigma_{10}} {\Sigma^\top _{10}\atop \Sigma_{11}} \right)
$$
 is a block matrix,
 $$\Sigma_{00}=\big(K(t_i,t_j)\big)_{i,j=1}^N,~
\Sigma_{10}=\Big(\frac{\partial}{\partial t_i} K(t_i,t_j)\Big)_{i,j=1}^N,~
\Sigma_{11}=\Big(\frac{\partial^2}{\partial t_i \partial t_j} K(t_i,t_j)\Big)_{i,j=1}^N
$$
are $N\!\times\! N$-matrices,
$$K_{t_0,{2N}}=\Big (K(t_0,t_1),\ldots,K(t_0,t_N),
\frac{\partial}{\partial t_0}K(t_0,t_1),\ldots,\frac{\partial}{\partial t_0}K(t_0,t_N)\Big )^\top $$
is a vector in $\mathbb{R}^{2N}$,
$X_{2N}=(f(t_1),\ldots,f(t_N),f'(t_1),\ldots,f'(t_N))^\top $ is an $2N\!\times\! m$-matrix and
$$
\hat\theta_{\mathrm{BLUE,}{2N}}=(X^\top _{2N}\Sigma^{-1}_dX_{2N})^{-1}X^\top _{2N}\Sigma^{-1}_dY_{2N}
$$
is the BLUE of $\theta$.
The MSE of the BLUP \eqref{eq:kriging-deriv} is given by
\bea
 \mathrm{MSE}(\hat y(t_0))=K(t_0,t_0)-\left[\begin{matrix}f(t_0)\\ K_{t_0,{2N}}\end{matrix}\right]^\top
 \left[\begin{matrix}0&X^\top _{2N}\\X_{2N}&\Sigma_{2N} \end{matrix}\right]^{-1}
 \left[\begin{matrix}f(t_0)\\ K_{t_0,{2N}}\end{matrix}\right].
\eea

For more general cases of prediction of processes and fields with derivatives observed at a finite number of points,  see \citep{morris1993bayesian,nather2003effective}.

\subsection{Continuous observations on an interval}
\label{sec32}

Consider  the  continuous-time model \eqref{1.1}, where
the error process $\epsilon$  has~a $q$ times differentiable covariance kernel $K(t,s)$. We also assume that  the vector-function $f$ is
 $q$ times differentiable and therefore  the response $y$ is $q$ times differentiable as well.

Suppose  we observe realization  $y(t)=y^{(0)}(t)$ for $t \in \mbox{\textsc{T}}_0 \subset \mathbb{R}$ and assume  that observations of the derivatives  $y^{(i)}(t)$ are also available
 for all $t\in\mbox{\textsc{T}}_i$, where  $\mbox{\textsc{T}}_i \subset \mathbb{R}$; $i=1, \ldots, q$. The sets $\mbox{\textsc{T}}_i$ ($i=0,1, \ldots, q$) do not have to be the same; some of these sets (but not all) can even be empty. If at least one of the sets $\mbox{\textsc{T}}_i$ contains an interval then we  speak of a problem with continuous observations.

Consider the problem of prediction  of $y^{(p)}(t_0)$, the $p$-th derivative of $y$ at a point $t_0\not\in\mbox{\textsc{T}}_p$, where $0 \leq p \leq q$.

A  general linear  predictor of the $p$-th derivative $y^{(p)}(t_0)$ can be defined as
\be
\label{eq:est2}
 \hat y_{p,Q}(t_0)= \int  \mathbf{Y}^\top (t)\mathbf{Q}(dt)
 = \sum_{i=0}^q \int_{\mbox{\textsc{T}}_i} y^{(i)}(t)Q_i(dt),
\ee
where $\mathbf{Y}(t)=\left( y(t), y^{(1)}(t), \ldots, y^{(q)}(t)\right)^\top $ is a vector with observations of the process and its derivatives,
$\mathbf{Q}(dt)= (Q_0(dt),\ldots,Q_q(dt))^\top $ is a vector of length $ (q+1)$ and
$Q_0(dt),\ldots,Q_q(dt)$ are signed measures
defined on~$\mbox{\textsc{T}}_0,\ldots, \mbox{\textsc{T}}_q$, respectively. The covariance matrix of $\mathbf{Y}(t)$ is
$$
\mathbb{K}(t,s)=E[ \mathbf{Y}(t) - E\mathbf{Y}(t)][ \mathbf{Y}(t) - E\mathbf{Y}(t)]^\top=\bigg(\frac{\partial^{i+j}K(t,s)}{\partial t^i\partial s^j}\bigg)_{i,j=0}^q
$$
which is a non-negative definite matrix of size $(q+1) \times (q+1) $.

The estimator  $\hat y_{p,Q}(t_0) $ is unbiased if $\mathbb{E}[\hat y_{p,Q}(t_0)]=\mathbb{E}[y^{(p)}(t_0)]$, which is equivalent to
$$\int  \mathbf{F}(t) \mathbf{Q}(dt)= f^{(p)}(t_0),$$
where  $\mathbf{F}(t)=\left( f(t), f^{(1)}(t), \ldots, f^{(q)}(t)\right)$ is a $m \!\times\! (q+1)$-matrix. \\


{\bf Assumption A$^{\prime\prime}$.} {\it

(1) The best linear unbiased estimator (BLUE)
  $\hat{\theta}_{\mathrm{BLUE}} = \int  \mathbf{G}(dt)\mathbf{Y}(t) $  exists in the model \eqref{1.1}, where $\mathbf{G}(dt)$ is some  signed $m \!\times\! (q+1)$-matrix measure (that is, the $j$-th column of $\mathbf{G}(dt)$ is a signed vector measure defined on ${\T}_j$);

 (2) There exists
 a signed vector-measure $\zeta_{p,t_0}(dt)$ (of size $q+1$) which satisfies  the equation
\be
\int   \mathbf{K}^\top(t,s)\zeta_{p,t_0}(dt)=
 \frac{\partial^{p}K(s,t_0)}{\partial t_0^p},\;\; \forall s\in \mbox{\textsc{T}}_i\, ,
\label{eq:zeta-p}
\ee
where $\mathbf{K}(t,s)=\big(\frac{\partial^{j}K(t,s)}{\partial s^j}\big)_{j=0}^q$ is a
$(q+1)$-dimensional  vector.
}

\vspace{2mm}
The problem of existence and construction of the BLUE in the continuous model with derivatives is discussed in \citep{DPZ2018}.
A general statement establishing the existence and explicit form of the BLUP is as follows. The proof is given in Section \ref{sec4}.

\begin{theorem}
\label{th:predict-t0-p}
If  Assumption A$^{\prime\prime}$ holds,
then
 the  BLUP measure $\mathbf{Q}_{*}$ exists and  is given by
\be
\label{eq:blup}
\mathbf{Q}_{*}(dt)=\zeta_{p,t_0}(dt)+\mathbf{G}^\top (dt) c_p,
\ee
where the  signed measure $\zeta_{p,t_0}(dt)$ satisfies \eqref{eq:zeta-p}
and
$$
c_p= f^{(p)}(t_0)-\int  \mathbf{F}(t)\zeta_{p,t_0}(dt) .
$$
The MSE  of the BLUP  $\hat y_{p,Q_*}(t_0)$ is given by
\bea
\mathrm{MSE}(\hat{y}_{p,Q_*}(t_0))=\left.\frac{\partial^{2p} K(t,s)}{\partial t^p \partial s^p}\right|_{{t=t_0 \atop s=t_0}}
+c_p^\top D f^{(p)}(t_0)- \int  \mathbf{K}^\top (t,t_0)\mathbf{Q}_{*}(dt)\, ,
\eea
where
$$
D= \int \!\! \int   \mathbf{G}(dt) \mathbb{K}(t,s) \mathbf{G}^\top (ds)
$$
is the covariance matrix of  $\hat{\theta}_{\mathrm{BLUE}}=\int  \mathbf{G}(dt)\mathbf{Y}(t) $.
\end{theorem}

\begin{example}  \label{ex3}
{ \rm  As a particular case of prediction in the model \eqref{1.1}, in this example we consider the problem  of predicting  a value of a process (so that $p=0$) with Mat\'{e}rn   $3/2$ covariance kernel
$
 K(t,s)=(1+\lambda|t-s|)e^{-\lambda|t-s|}\,;
$
this kernel is once differentiable and is very popular in practice, see e.g. \citep{rasmussen2006gaussian}. We assume that the vector-function $f$ in the model \eqref{1.1} is 4 times differentiable and that the process $y$ and its derivative $y^\prime$ are observed on an interval $[A,B]$ (so that $\mbox{\textsc{T}}_0=\mbox{\textsc{T}}_1=[A,B]$ in the general statements).
As shown in \citep{DPZ2018}, for this kernel  the BLUE measure $\mathbf{G}(dt)$ can be expressed in terms of the signed matrix-measure $\zeta(dt)=(\zeta_0(dt),\zeta_1(dt))$ with
 \begin{eqnarray*}
 \zeta_0(dt)& =&z_A\delta_A(dt)+z_B\delta_B(dt)+z(t)dt ,  \\
\zeta_1(dt) &=&z_{1,A}\delta_A(dt)+z_{1,B}\delta_B(dt),
\end{eqnarray*}
where
\bea
 z_A&=&\frac{1}{4\lambda^3}\big(f^{(3)}(A) -3\lambda^2 f^{(1)}(A)  + 2\lambda^3 f(A)\big),\nonumber\\
 z_{1,A}&=&\frac{1}{4\lambda^3}\big(-f^{(2)}(A)+ 2\lambda f^{(1)}(A)- \lambda^2 f(A)\big),\nonumber\\
  z_B&=&\frac{1}{4\lambda^3}\big(-f^{(3)}(B) +3\lambda^2 f^{(1)}(B)  + 2\lambda^3 f(B)\big),\\
  z_{1,B}&=&\frac{1}{4\lambda^3}\big(f^{(2)}(B)+ 2\lambda f^{(1)}(B)+ \lambda^2 f(B)\big),\nonumber\\
  z(t)&=&\frac{1}{4\lambda^3 }\big( \lambda^4 f(t)-2\lambda^2f^{(2)}(t)+f^{(4)}(t)\big).\nonumber
\eea
Then using \citep[Sect. 3.4]{DPZ2018} we obtain $\zeta_{0,t_0}(dt)=(\zeta_{0,t_0,0}(dt),$ $ \zeta_{0,t_0,1}(dt))$ with
\bea
\zeta_{0,t_0,0}(dt)&=&z_{t_0,A}\delta_A(dt)+z_{t_0,B}\delta_B(dt)+z_{t_0}(t)dt,\\
\zeta_{0,t_0,1}(dt)&=&z_{t_0,1,A}\delta_A(dt)+z_{t_0,1,B}\delta_B(dt),
\eea
where for  $t_0>B$ we have $z_{t_0,A}=0$, $z_{t_0,1,A}=0,$ $ z_{t_0}(t)=0,$
\bea
  z_{t_0,B}=(1+\lambda(t_0-B)) e^{-\lambda(t_0-B)},\;\;\;
   z_{t_0,1,B}=(t_0-B) e^{-\lambda(t_0-B)}\, .
\eea

We also obtain the matrix $$C=\int_A^B \zeta_0(dt)f^\top(t)+\int_A^B \zeta_1(dt){f'}^\top(t)$$
defined in \citep[Lem. 2.1]{DPZ2018} from the condition of unbiasedness.  If $D$, the covariance matrix of the BLUE is non-degenerate, then  $D=C^{-1}$. In the present case,
 \bea
 C&=&\frac{1}{2}\Big[f(A)f^\top (A)+f(B)f^\top (B)\Big]+\frac{1}{2\lambda^2}\Big[f'(A)f'^\top (A)+f'(B)f'^\top (B)\Big]+\\
 &&-\frac{1}{4\lambda}\Big[f'(A)f^\top (A)+f(A)f'^\top (A)+f'(B)f^\top (B)+f(B)f'^\top (B)\Big]+\\
 &&+\frac{1}{4\lambda^3}\int_A^B \Big[\lambda^4f(t)f^\top (t)+2\lambda^2f'(t)f'^\top (t)+f''(t)f''^\top (t)\Big]dt\, ,\\
 &&c_0=\Big(f(t_0)-[z_{t_0,B}f(B)+z_{t_0,1,B}f'(B)]\Big).
\eea
The BLUE-defining measure $\mathbf{G}(dt)$ is expressed through the measures $\zeta(dt)$ and  the matrix $C$ by $\mathbf{G}(dt)= C^{-1}\zeta(dt)$.
The BLUP measure for process prediction is given by
\bea
\mathbf{Q}_{*}(dt)&=&\zeta_{0,t_0}(dt)+\mathbf{G}^\top (dt) c_0\\
&=&\Big(\zeta_{0,t_0,0}(dt)+c_0^\top C^{-1}\zeta_0(dt),\zeta_{0,t_0,1}(dt)+c_0^\top C^{-1}\zeta_1(dt)\Big)^\top,
\eea
where
$$
c_0= f(t_0)-\int  \mathbf{F}(t)\zeta_{0,t_0}(dt) .
$$

For the location scale model with $f(t)=1$, we obtain $C=1+\lambda (B-A)/4$, $c_0=(1-z_{t_0,B})$ and, therefore, a BLUP measure for this model is given by
\be
 \nonumber
\mathbf{Q}_{*}(dt)&=& 0.5c_0\delta_A(dt)/C+(0.5c_0/C+z_{t_0,B})\delta_B(dt)+0.25c_0\lambda dt/C\\
&& ~~-0.25c_0/(C\lambda)\delta_A(dt) +(z_{t_0,1,B}+0.25c_0/(C\lambda))\delta_B(dt)
 \nonumber
\ee
Therefore, the corresponding BLUP  is given
\bea
\hat y_{0,Q_*}(t_0)&=&0.5c_0y(A)/C+(0.5c_0/C+z_{t_0,B})y(B)+0.25c_0\lambda \int_A^B y(t)dt/C\\
&& -0.25c_0/(C\lambda) y'(A)(dt) +(z_{t_0,1,B}+0.25c_0/(C\lambda))y'(B).
\eea

Table~\ref{tab:tt} gives values of the square root of the MSE of the BLUP in the location scale model
at the point $t_0=2$ for three families of designs, where $[A,B] = [0,1]$.
We observe that
observations of derivatives inside the interval do not bring any improvement to the BLUP which can be explained by the fact that the  weights of the continuous BLUP at  derivatives at  points in the interior of the interval $[A,B]$ are $0$.

\begin{table}[!hhh]
\caption{\it The square root of the MSE of the BLUP at the point $t_0=2$ for different designs.
(i) the design $\xi_{N,0}$ observing the process at $N$-point equidistant points,
(ii) the design $\xi_{N,2}$ observing the process at $N$-point equidistant points
and the derivative at two boundary points.
(iii) the design $\xi_{N,N}$ observing the process and derivative at $N$-point equidistant points. The model
is  the location scale model on the interval $[0,1]$
and the covariance kernel of the error process is given by  the
Mat\'{e}rn   3/2 kernel with $\lambda=2$. For continuous observations the square root of the BLUB is given by $\sqrt{\mathrm{MSE}}=0.9985569896$.}
\begin{center}
\begin{tabular}{|c|c|c|c|c|c|}
\hline
 $N$ & 2 & 4 & 8 & 16 \\
\hline
$\xi_{N,0}$& 1.059339& 1.038152& 1.019244& 1.009052  \\
$\xi_{N,2}$&0.999276&0.9985675343&0.9985573516&0.9985570068 \\
 $\xi_{N,N}$&0.999276&0.9985675343&0.9985573516&0.9985570068 \\
\hline
\end{tabular}
\end{center}
\label{tab:tt}
\end{table}
}
 \end{example}

\subsection{Location scale model on a product set}
\label{sec:sec33}

Similarly to Section~\ref{sec24}, we consider the location scale model
\eqref{eq:loc_model}
defined on the product set $\T=\T_1 \times \T_2$ (where $\T_1$ and $\T_2$ are Borel sets in $\mathbb{R}$)
with the kernel $ \mbox{\textsf{K}}$ of the random field ${\ve}$  having the product form
\eqref{eq:loc_prod}.
The results of this section (as of Section~\ref{sec24}) can be easily generalized to the case of $d>2$ variables.

Assume that Assumption A$^{\prime\prime}$ with $q=1$ is satisfied for two one-dimensional models~\eqref{eq:prod_models}. For this assumption to hold, the  process  $ \{  y(t_1,t_2) ~|~  (t_1,t_2) \in \T \} $
 has to be once differentiable with respect to $t_1$ and $t_2$.
 Let the  measures $G_{0,i}(du)$ and $G_{1,i}(du)$
define  the BLUE
$$
\int_{\T_i} y_{(i)}(u) G_{0,k}(du)+\int_{\T_i} y_{(i)}^{(1)}(u) G_{1,i}(du)
$$
in the univariate models \eqref{eq:prod_models}; $i=1,2$.
 In this case, results of \citep{DPZ2018} imply that the BLUE of $\theta$ in the model \eqref{eq:loc_model} has the form
$
\hat {\theta}=
\int_\T  {\mathbf{Y}}^\top (t) \mbox{ {\bf  \textsf{G}}}(dt),
$
where
$${\mathbf{Y}}^\top (t)=\left(  y(t),\frac{\partial}{\partial t_1}  y(t),
\frac{\partial}{\partial t_2}  y(t), \frac{\partial^2}{\partial t_1\partial t_2} y(t)\right)$$
and
$$
\mbox{ {\bf  \textsf{G}}}(dt)=\left( \mbox{\textsf{G}}_{00}(dt),\mbox{\textsf{G}}_{10}(dt),\mbox{\textsf{G}}_{01}(dt),\mbox{\textsf{G}}_{11}(dt)\right)^\top
$$
with $\mbox{\textsf{G}}_{ij}(dt)=G_{i,k}(dt_1)G_{j,k}(dt_2)$.

Assume we want to predict $ y(T)$ at a point $T=(T_1,T_2) \notin \T$.
The analogue of the equation \eqref{eq:zeta1} is given by
\be
\label{eq:eqnew}
 \int_\T  \mbox{\bf\textsf{K}}^\top (t,t') \mbox{\bf\textsf{Z}}_{T}(dt')=  \mbox{\textsf{K}}(t,T),\;\; \forall t\in \T,
\ee
where
\bea
 \mbox{\bf\textsf{K}}((t_1,t_2),(s_1,s_2))=
 \left(                    \begin{array}{l}
 K_1(t_1,s_1)K_2(t_2,s_2)\\
 \frac{\partial}{\partial t_1}K_1(t_1,s_1)K_2(t_2,s_2)\\
 K_1(t_1,s_1)\frac{\partial}{\partial t_2}K_2(t_2,s_2)\\
 \frac{\partial}{\partial t_1}K_1(t_1,s_1)\frac{\partial}{\partial t_2}K_2(t_2,s_2)
\end{array}\right).
\eea

Observing the product-form of expressions, we directly obtain that
a solution of \eqref{eq:eqnew}  has the form
\bea
\mbox{\bf\textsf{Z}}_{T}(dt_1,dt_2)=
\left(                    \begin{array}{l}
\zeta_{0,T_1}(dt_1) \zeta_{0,T_2}(dt_2)\\ \zeta_{1,T_1}(dt_1) \zeta_{0,T_2}(dt_2)\\
\zeta_{0,T_1}(dt_1) \zeta_{1,T_2}(dt_2)\\ \zeta_{1,T_1}(dt_1) \zeta_{1,T_2}(dt_2)
\end{array}\right),
\eea
where measures $\zeta_{0,T_i}(dt)$ and $\zeta_{1,T_i}(dt)$ for $i=1,2$ satisfy the equation
\bea
 \int_{\T_i}   K_i(t,s) \zeta_{0,T_i}(dt)+\int_{\T_i}   \frac{\partial}{\partial t}K_i(t,s) \zeta_{1,T_i}(dt)=  K_i(s,T_i) ,\;\; \forall s\in {\T_i}.~
\eea

Finally, the BLUP at the point $T=(T_1,T_2)$ is
$\int_\T {\mathbf{Y}}^\top (t) \mbox{\bf\textsf{Q}}_{*}(dt) $,
where
$
 \mbox{\bf\textsf{Q}}_{*}(dt)= \mbox{\bf\textsf{Z}}_{T}(dt)+ {c}_0 \mbox{{\bf\textsf{G}}}(dt)
$
with $
c_0=1- \int_\T (1,0,0,0)\mbox{\bf\textsf{Z}}_{T}(dt).
$

The MSE of the BLUP is given by
$$
\mathrm{MSE}(\hat{ y}_{ 0, \mbox{\footnotesize\bf\textsf{Q}}_{*}  } (T))
=1+c_0D-\int_\T \mbox{\bf\textsf{K}}^\top (t,T)
\mbox{\bf\textsf{Q}}_{*}(dt),
$$
where $D$ is the variance of  the BLUE.

\begin{example}
{\rm Consider a   location scale model on a square  $[0,1]^{2}$
with a product covariance Mat\'{e}rn   $3/2$  kernel, that is
\bea
 \mbox{\textsf{K}}(t,t')&=&\mathbb{E}[\ve(t)\ve(t')]=K(t_1,t_1')K(t_2,t_2')
,
\eea
where
\begin{equation} \label{mat}
 K(u,u')=(1+\lambda|u-u'|)e^{-\lambda|u-u'|}.
\end{equation}
Define the measures
$$G_0(du) = \frac1{4+\lambda} \left[  2\delta_0 (du)+ 2\delta_1 (du)+\lambda  du \right]$$
and
$$
G_1(du)=\frac1{(4+\lambda)\lambda} \left[  \delta_1 (du)- \delta_0 (du)\right],\;\;u \in [0,1].
$$
In view of \citep[Sect. 3.4]{DPZ2018},
$$
\int_0^1 y(u) G_0(du)+\int_0^1 y^{(1)}(u) G_1(du)
$$
defines a   BLUE in the model $y(u)=\theta + \ve(u)$ with $u \in [0,1]$  and  covariance kernel
\eqref{mat}.
Additionally, from \citep[Sect. 3.4]{DPZ2018} we have
\bea
 \zeta_{0,T_i}(du)= \left\{
                    \begin{array}{ll}
                     (1+\lambda|T_i|)e^{-\lambda |T_i|} \delta_0(du),  & T_i \leq 0, \\
                      \delta_{T_i}(du),   & 0 \leq T_i \leq 1,  \\
                      (1+\lambda(T_i-1))e^{-\lambda (T_i-1)} \delta_1(du), & T_i \geq 1, \\
                    \end{array}
                  \right.
\eea
and
\bea
 \zeta_{1,T_i}(du)= \left\{
                    \begin{array}{ll}
                     -|T_i|e^{-\lambda |T_i|} \delta_0(du),  & T_i \leq 0, \\
                      0,   & 0 \leq T_i \leq 1,  \\
                      (T_i-1)e^{-\lambda (T_i-1)} \delta_1(du), & T_i \geq 1. \\
                    \end{array}
                  \right.
\eea

Finally, $c_0=1- \int_0^1\!\! \int_0^1 (1,0,0,0)\mbox{\bf\textsf{Z}}_{T}(dt)=
1-\int_0^1\!\zeta_{0,T_1}(dt_1) \! \int_0^1  \zeta_{0,T_2}(dt_2)$ and the BLUP measure is given by
$
\mbox{\bf\textsf{Q}}_{*}(dt)= \mbox{\bf\textsf{Z}}_{T}(dt)+ c_0 \mbox{{\bf\textsf{G}}}(dt);
$ that is,
\bea
\mbox{\bf\textsf{Q}}_{*}(dt)=
\left(                    \begin{array}{l}
\zeta_{0,T_1}(dt_1) \zeta_{0,T_2}(dt_2)+c_0 G_0(dt_1)G_0(dt_2)\\
\zeta_{1,T_1}(dt_1) \zeta_{0,T_2}(dt_2)+c_0 G_1(dt_1)G_0(dt_2)\\
\zeta_{0,T_1}(dt_1) \zeta_{1,T_2}(dt_2)+c_0 G_0(dt_1)G_1(dt_2)\\
\zeta_{1,T_1}(dt_1) \zeta_{1,T_2}(dt_2)+c_0 G_1(dt_1)G_1(dt_2)
\end{array}\right).
\eea
We now investigate the performance of five discrete designs:
\begin{itemize}
  \item[(i)]
 the design $\xi_{N^2,0,0,0}$, where we observe process $y$  on an $N\!\times\! N$ grid;
  \item[(ii)]
the design $\xi_{N^2,4,4,4}$, where we observe process $y$  on an $N\!\times\! N$ grid
 and additionally derivatives
$\frac{\partial y}{\partial t_1}$, $\frac{\partial y}{\partial t_2}$,
$\frac{\partial^2 y}{\partial t_1\partial t_2}$ at 4 corners of $[0,1]^2$;
  \item[(iii)]
the design $\xi_{N^2,N^2,N^2,0}$, where we observe process $y$ and derivatives
$\frac{\partial y}{\partial t_1}$, $\frac{\partial y}{\partial t_2}$
 on an $N\!\times\! N$ grid;
  \item[(iv)]
  the design $\xi_{N^2,N^2,N^2,0}$,
  where we observe process $y$  on an $N\!\times\! N$ grid
   and derivatives
$\frac{\partial y}{\partial t_1}$, $\frac{\partial y}{\partial t_2}$, $\frac{\partial^2 y}{\partial t_1\partial t_2}$
at $4N-4$ equidistant points on the boundary of $[0,1]^2$;
  \item[(v)]
the design $\xi_{N^2,N^2,N^2,N^2}$, where we observe process $y$ and derivatives
$\frac{\partial y}{\partial t_1}$, $\frac{\partial y}{\partial t_2}$,
$\frac{\partial^2 y}{\partial t_1\partial t_2}$ at $N\!\times\! N$ equidistant points on an $N\!\times\! N$ grid.
\end{itemize}

\begin{table}[!hhh]
\caption{\it The square root of the MSE of the BLUP at the point $(2,2)$  (upper part)
and the point  $(0.5,2)$  (lower part)  for several designs
in the location scale model on the square $[0,1]^2$
with Mat\'{e}rn   3/2 product-kernel ($\lambda=2$).
The square root of the MSE of the continuous BLUP equals 1.119510 at the point $(2,2)$
and 0.958494 at the point  $(0.5,2)$.}
\begin{center}
\begin{tabular}{|l|c|c|c|c|c|c|}
\hline
 $\;\;\;\;\;\;\;\;\;\;\;\;\; N$ & 2 & 3 & 4 & 8 & 16  \\
\hline
 $\xi_{N^2,0,0,0}$         &1.16139& 1.15344 & 1.14972 &1.13548& 1.12764 \\
 $\xi_{N^2,4,4,4}$         &1.121205&1.119682&1.119582&1.119543&1.119528\\
 $\xi_{N^2,N^2,N^2,0}$     &1.124401&1.121576&1.120913&1.119893&1.119609\\
 $\xi_{N^2,4N-4,4N-4,4N-4}$&1.121205&1.119632&1.119535&1.119511&1.119510\\
\hline\hline
$\xi_{N^2,0,0,0}$&1.03152& 1.00413 & 0.99900 &0.97862& 0.96862 \\
 $\xi_{N^2,4,4,4}$         &0.979953&0.962754&0.963426&0.960604&0.959550\\
 $\xi_{N^2,N^2,N^2,0}$     &0.982184&0.958732&0.959663&0.958606&0.958511\\
 $\xi_{N^2,4N-4,4N-4,4N-4}$&0.979953&0.958566&0.959314&0.958556&0.958500\\
\hline
\end{tabular}
\end{center}
\label{tab3}
\end{table}

The results  are depicted in Table \ref{tab3}, which shows the square root of the MSE of predictions at the point $(2,2)$
and  $(0.5,2)$ for different sample sizes.
 For any given $N\geq 2$, the MSE for prediction outside the square $[0,1]^2$ for the designs $\xi_{N^2,N^2,N^2,N^2}$ and $\xi_{N^2,4N-4,4N-4,4N-4}$ are exactly the same. This is related to the fact  that  the BLUP weights associated with all derivatives at interior points in $[0,1]^2$ of the designs $\xi_{N^2,N^2,N^2,N^2}$ are all 0. This means that for optimal prediction of $y(t_0)$ at a point $t_0$ outside the observation region  one needs the design guaranteeing the optimal BLUE plus the observations of $y(t)$ and $y'(t)$ at points $t$ closest to $t_0$. Note that the results of \citep[Sect. 3.4]{DPZ2018} imply that the continuous  optimal design for the BLUE does not use values of any derivatives of the process (or field for the product-covariance model) in the interior of $\T$.

 The observation above  is consistent with our other numerical experience which have shown that the BLUP at a point $t_0 \in (0,1) \times (0,1)$ constructed from the  design $\xi_{N^2,N^2,N^2,N^2}$
has vanishing weights at  all derivatives of interior points of  $[0,1]^2$ with five exceptions:
the center $0$  and the four points which are closest to $t_0$ in the $L_\infty$ (Manhattan) metric.

\begin{figure}[!hhh]
\begin{center}
 \includegraphics[width=0.48\textwidth]{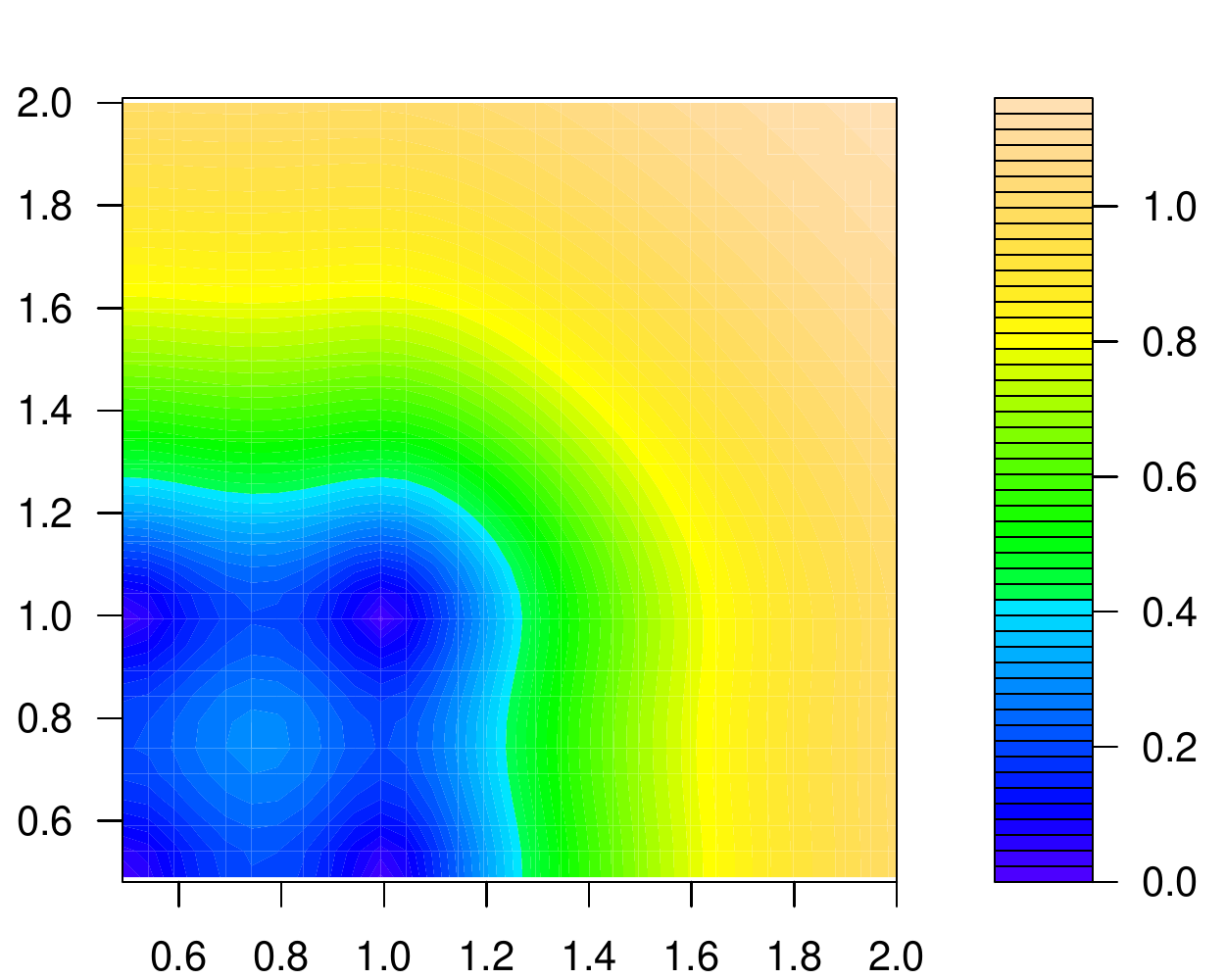}
 \includegraphics[width=0.48\textwidth]{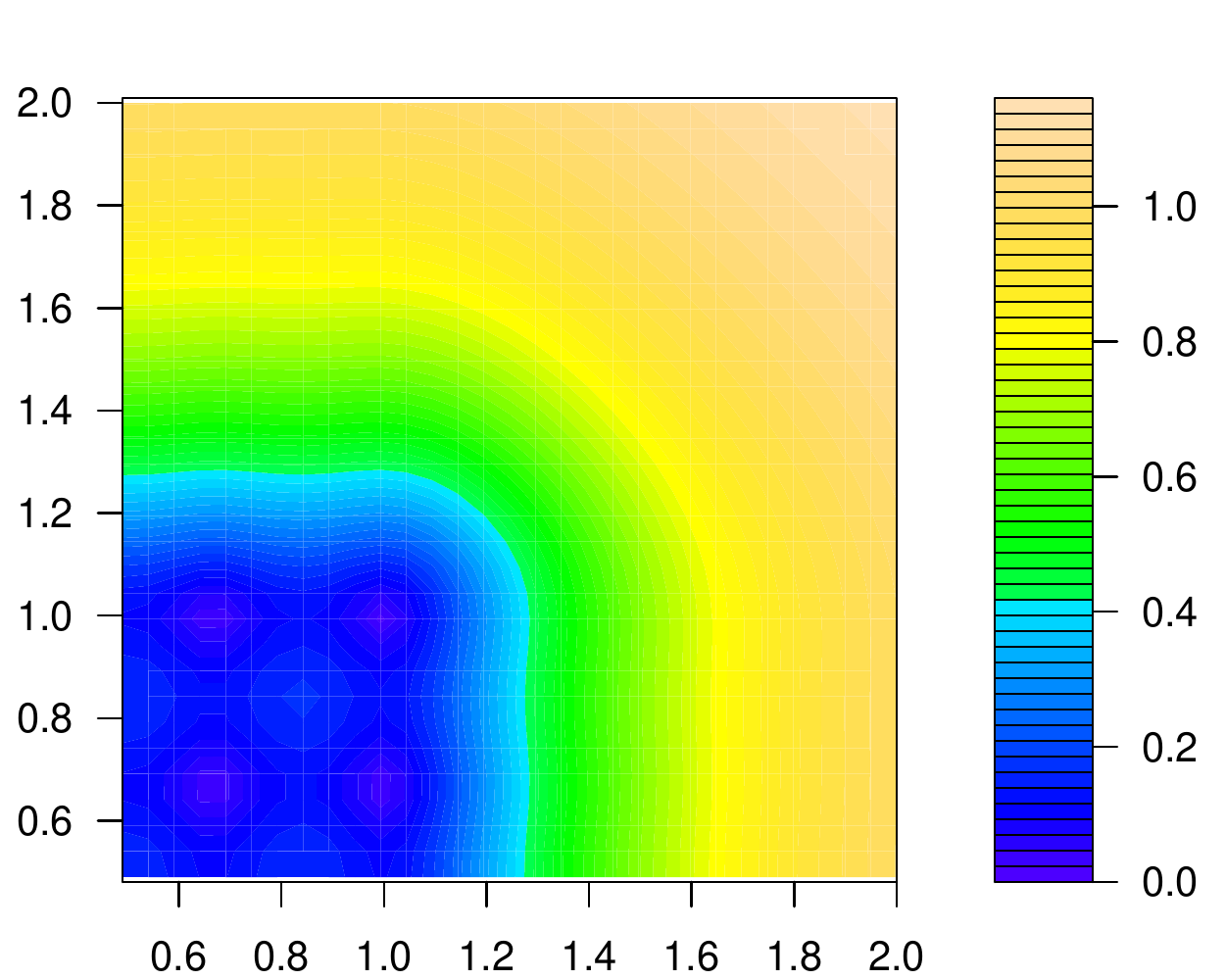}
\end{center}
\caption{\it  Square root of the MSE of the BLUP for the design $\xi_{N^2,0,0,0}$
with $N=3$ (left) and $N=4$ (right), and the Mat\'{e}rn   3/2 product-kernel with $\lambda=2$.}
\label{fig:rmse-matern}
\end{figure}

\begin{figure}[!hhh]
 \includegraphics[width=0.48\textwidth]{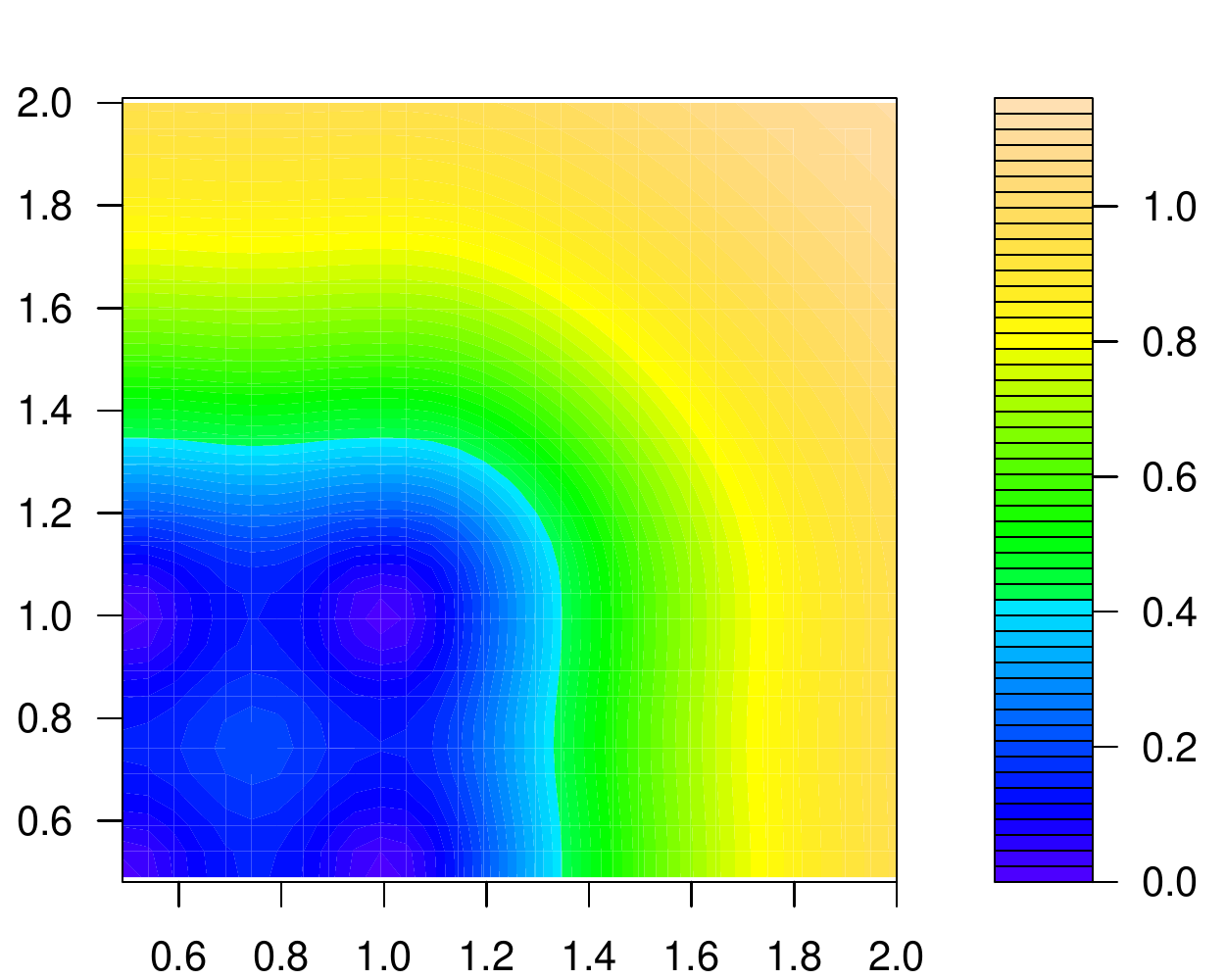}
 \includegraphics[width=0.48\textwidth]{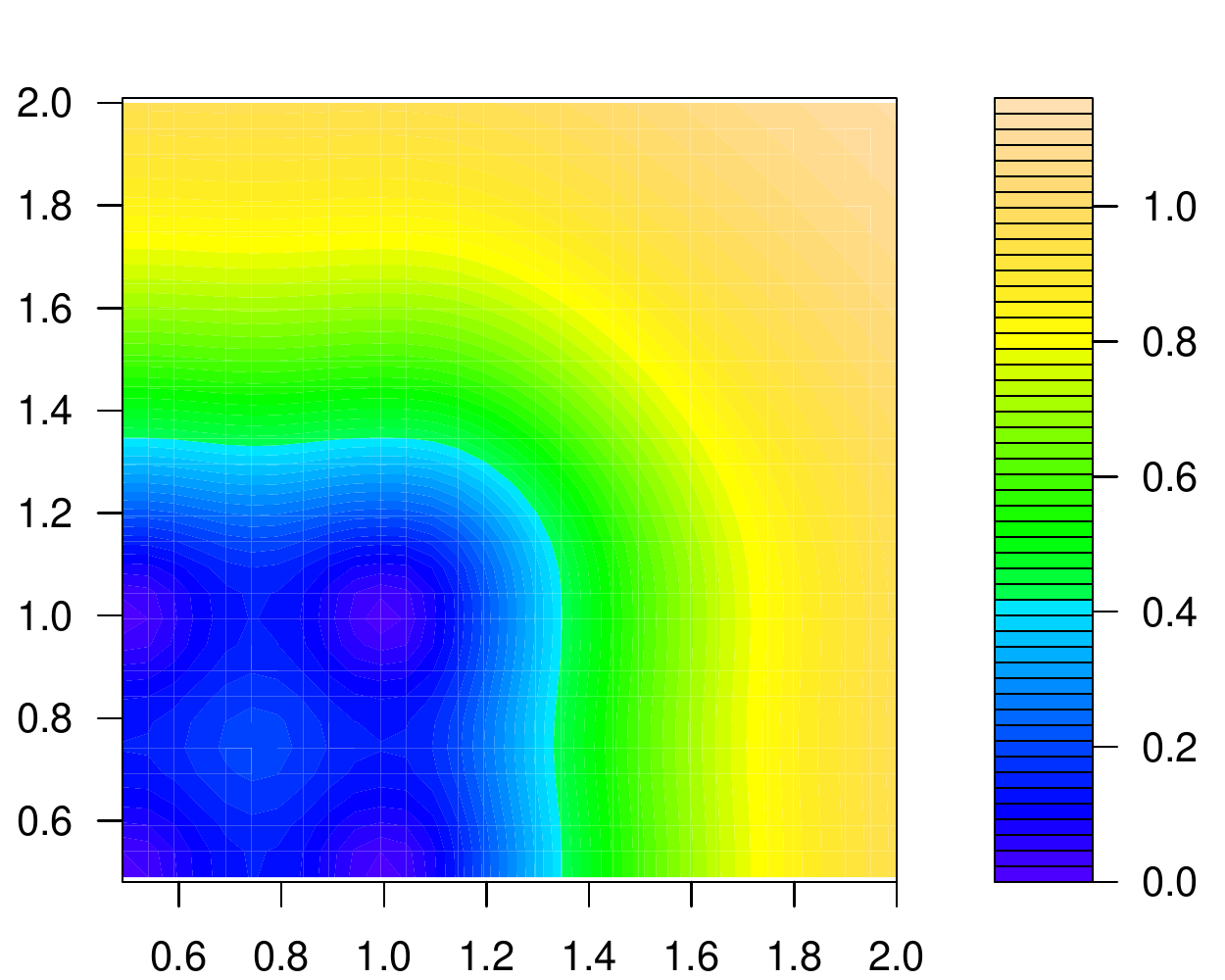}
\caption{\it  Square root of the MSE of the BLUP for the design $\xi_{N^2,4N-4,4N-4,4N-4}$  (left) and
$\xi_{N^2,N^2,N^2,N^2}$ (right) with $N=3$ and the Mat\'{e}rn   3/2 product-kernel with $\lambda=2$.}
\label{fig:rmse-matern2}
\end{figure}

Figures \ref{fig:rmse-matern} and \ref{fig:rmse-matern2} compare the MSE for some designs.
As Figure \ref{fig:rmse-matern2} illustrates, additionally to Table \ref{tab3}, the MSE for designs
$\xi_{N^2,4N-4,4N-4,4N-4}$   and
$\xi_{N^2,N^2,N^2,N^2}$
is exactly the same for all points outside $[0,1]^2$ and
almost the same at all interior points of $[0,1]^2$.

}
\end{example}


\section{Proofs} \label{sec4}
\def\theequation{4.\arabic{equation}}
\setcounter{equation}{0}

\subsection{Proof of Theorem \ref{th:predict-nu}}

To start, we proof the following lemma.

\begin{lemma}
The mean squared error [relative to the true process value] of any unbiased estimator $\hat Z_Q=  \int_\T y(t)Q(dt)$  is given by
\bea
 &&\mathrm{MSE}(\hat Z_{Q})=\mathbb{E}\left(Z -\hat Z_{Q}\right)^2 \label{Feq:var}=\\
&&\!\!\! \int_{\mathcal{S}}\!\int_{\mathcal{S}}\!K(t,s)\nu(dt)\nu(ds)\!-\!2\!\int_{\mathcal{S}}\!\int_\T\! K(t,s)\nu(dt)Q(ds)
 \! +\!\int_\T \!\int_\T\! Q(dt)K(t,s)Q(ds)  .\nonumber
\eea
\end{lemma}

{\bf Proof.} Straightforward calculation gives
\bea
 \!\!\mathrm{MSE}(\hat Z_Q)\!\!&=&\mathbb{E}\left(Z -\hat Z_{Q}\right)^2
 =\mathbb{E}\left(Z- \int_\T y^{}(t)Q(dt) \right)^2\\
 &=&\mathbb{E}\left(\int_{\mathcal{S}} [\theta^\top f(t)+\epsilon(t)]\nu(dt) - \int_\T [\theta^\top f^{}(t)+\epsilon^{}(t)]Q(dt) \right)^2\\
 &=&\mathbb{E}\left(\int_{\mathcal{S}} \epsilon(t)\nu(dt) - \int_\T \epsilon^{}(t)Q(dt) \right)^2\\
 &=&\!\mathbb{E}\left(\int_{\mathcal{S}} \epsilon(t)\nu(dt)\! -\! \int_\T \epsilon^{}(t)Q(dt) \right)  \left(\int_{\mathcal{S}} \epsilon(s)\nu(ds) \!- \!\int_\T \epsilon^{}(s)Q(ds) \right)\\
  &=&\!\int_{\mathcal{S}}\int_{\mathcal{S}}K(t,s)\nu(dt)\nu(ds)\!-\!2\int_{\mathcal{S}}\int_\T K(t,s)\nu(dt)Q(ds) \\
  &&+ \int_\T \int_\T K(t,s)Q(dt)Q(ds)\, ,
\eea
as  required.
\hfill$\Box$

\medskip

Let us now prove the main result.
We will show that $\mathrm{MSE}(\hat Z_Q)\ge\mathrm{MSE}(\hat Z_{Q_*})$,
where $\hat Z_Q$ is any linear unbiased estimator of the from \eqref{zqaht}  and $\hat Z_{Q_*}$
is defined by the measure \eqref{eq:blup1}.
Define $R(dt)=Q(dt)-Q_{*}(dt)$. From the condition of unbiasedness for $Q(dt)$ and $Q_{*}(dt)$,
we have $\int_\T f(t)R(dt)=0_{m \times 1}$.

We obtain
\bea
 \!\!\!\mathrm{MSE}(\hat Z_Q)\!\!\!&=&\mathrm{MSE}(\hat Z_{Q_{*}+R}) \\
 &=&\int_{\mathcal{S}}\int_{\mathcal{S}}K(t,s)\nu(dt)\nu(ds) -2\int_{\mathcal{S}} \int_\T K(t,s)\nu(dt) [Q_{*}+R](ds)\\
   &&+\int_\T \int_\T [Q_{*}+R](dt){K}(t,s)[Q_{*}+R](ds)  \\
 &=&\!\!\mathrm{MSE}(Q_{*})\!-\!2\int_{\mathcal{S}}\! \int_\T\! K(t,s)\nu(dt)R(ds)\!+\!\int_\T \!\int_\T \!R(dt){K}(t,s)R(ds)\\
    &&+2\int_\T \int_\T Q_{*}(dt){K}(t,s)R(ds)\\
 &\ge&\!\!\mathrm{MSE}(Q_{*})\!-\!2\!\int_{\mathcal{S}}\!\int_\T\! K(t,s)\nu(dt)R(ds)\!+\!2\!\int_\T \!\int_\T\! Q_{*}(dt){K}(t,s)R(ds)\\
 &=&\mathrm{MSE}(Q_{*})+2\int_\T \Big[\int_\T Q_{*}(dt){K}(t,s)- \int_{\mathcal{S}} K(t,s)\nu(dt)\Big]R(ds)\\
 &=&\mathrm{MSE}(Q_{*})+2\int_\T c^\top D f(s)R(ds) =\mathrm{MSE}(Q_{*})\, ,
\eea
where the inequality follows from nonnegative definiteness of the covariance kernel and the last equality follows from the unbiasedness condition
$\int f(t)R(dt)=0$.
\hfill$\Box$

\subsection{Proof of Theorem \ref{th:predict-t0-p}}

For simplicity, assume $p=0$; the case $p>0$ can be dealt with analogously.
First, we derive the following lemma.

\begin{lemma}
The mean squared error of any unbiased estimator $\hat y_{Q}(t_0)$ of the form  \eqref{eq:est2} is given by
\bea
 \!\!\mathrm{MSE}(\hat y_{Q}(t_0))\!\!&=&\mathbb{E}\left(y(t_0) -\hat y_{Q}(t_0)\right)^2 \label{eq:var}\\
  &=&
 K(t_0,t_0)-2\int_\T \mathbf{K}^\top (t_0,s)\mathbf{Q}(ds)+
 \int_\T \int_\T \mathbf{Q}^\top (dt)\mathbb{K}(t,s)\mathbf{Q}(ds) \, .\nonumber
\eea
\end{lemma}

{\bf Proof.} Straightforward calculation gives

\bea
 \!\!\!\!\mathrm{MSE}(\hat y_Q(t_0))\!\!\!&=&\mathbb{E}\left(y(t_0) -\hat y_{Q}(t_0)\right)^2 =\mathbb{E}\left(y(t_0)-\sum_{i=0}^q \int_\T y^{(i)}(t)Q_i(dt) \right)^2\\
 &=&\mathbb{E}\left(\theta^\top f(t_0)+\epsilon(t_0) -\sum_{i=0}^q \int_\T [\theta^\top f^{(i)}(t)+\epsilon^{(i)}(t)]Q_i(dt) \right)^2\\
 &=&\mathbb{E}\left(\epsilon(t_0) -\sum_{i=0}^q \int_\T \epsilon^{(i)}(t)Q_i(dt) \right)^2\\
  &=&K(t_0,t_0)-2\sum_{j=0}^q \int_\T \frac{\partial^{j}K(t_0,s)}{\partial s^j}Q_j(ds) \\
 &&+\sum_{i=0}^q \sum_{j=0}^q \int_\T \int_\T \frac{\partial^{i+j}K(t,s)}{\partial t^i\partial s^j}Q_i(dt)Q_j(ds)\, ,
\eea
as required.
\hfill$\Box$

Now we will prove the main result.
We will show that $\mathrm{MSE}(\hat y_Q(t_0))\ge\mathrm{MSE}(\hat y_{Q_{*}}(t_0))$,
where $\hat y_Q(t_{0})$ is any linear unbiased estimator of the form \eqref{eq:est2} and  $\hat y_{Q_{*}}(t_{0})$
is defined by \eqref{eq:blup}.
Define $\mathbf{R}(dt)=\mathbf{Q}(dt)-\mathbf{Q}_{*}(dt)$. From the condition of unbiasedness 
for $\mathbf{Q}(dt)$ and $\mathbf{Q}_{*}(dt)$,
we have
$
\label{eq:unbiasd}
\int_\T \mathbf{F}(t)\mathbf{R}(dt)=0_{m \times 1},
$
where $\mathbf{F}(t)=(f(t),f^{(1)}(t),\ldots,f^{(q)}(t))$. Therefore we  obtain
\bea
\!\!\! \mathrm{MSE}(\hat y_Q(t_0))\!\!\!&=&\mathrm{MSE}(\hat y_{Q_{*}+R}(t_0)) \\
 &=&K(t_0,t_0)-2\int_\T \mathbf{K}^\top (t_0,s)[\mathbf{Q}_{*}+\mathbf{R}](ds)\\
   &&+\int_\T \int_\T [\mathbf{Q}_{*}+\mathbf{R}]^\top (dt)\mathbb{K}(t,s)[\mathbf{Q}_{*}+\mathbf{R}](ds)  \\
 &=&\!\!\mathrm{MSE}(\hat y_{Q_{*}}(t_0))\!-\!2\!\int_\T\! \mathbf{K}^\top (t_0,s)\mathbf{R}(ds)\!+\!\int_\T\! \int_\T \!\mathbf{R}^\top (dt)\mathbb{K}(t,s)\mathbf{R}(ds)\\
    &&+2\int_\T \int_\T \mathbf{Q}_{*}^\top (dt)\mathbb{K}(t,s)\mathbf{R}(ds)\\
 &\ge&\!\!\mathrm{MSE}(\hat y_{Q_{*}}(t_0))\!-\!2\!\int_\T\! \mathbf{K}^\top (t_0,s)\mathbf{R}(ds)
 \!+\!2\!\int_\T \! \int_\T\! \mathbf{Q}_{*}^\top (dt)\mathbb{K}(t,s)\mathbf{R}(ds)\\
 &=&\mathrm{MSE}(\hat y_{Q_{*}}(t_0))+2\int_\T \Big [\int_\T \mathbf{Q}_{*}^\top (dt)\mathbb{K}(t,s)- \mathbf{K}^\top (t_0,s)\Big]\mathbf{R}(ds)\\
 &=&\mathrm{MSE}(\hat y_{Q_{*}}(t_0))+2\int_\T c_p^\top D\mathbf{F}(s)\mathbf{R}(ds)=\mathrm{MSE}(\hat y_{Q_{*}}(t_0))\, ,
\eea
where the inequality follows from nonnegative definiteness of the covariance kernel and the last equality follows from
the unbiasedness condition $\int \mathbf{F}(t)\mathbf{R}(dt)=0$.
\hfill$\Box$

\section{Appendix: more examples of predicting process values}
\def\theequation{5.\arabic{equation}}
\setcounter{equation}{0}

In the appendix, we give further examples of  prediction of values of specific random processes $y(t)$, which follows the model \eqref{1.1} and observed for all $t\in \T=[A,B]$.
In Section~\ref{sec:mark}, we illustrate application of  Proposition \ref{prop1} and in Section~\ref{sec:intB}
we give an example of application of Theorem~\ref{th:predict-t0-p}.  In the example of Section~\ref{sec:intB} we consider the integrated Brownian motion process, which is a once differentiable random process, and we assume that in addition to values of $y(t)$, the values of the derivative of $y(t)$ are also available. As in the main body of the paper, the components of the vector-function $f(t)$ in \eqref{1.1} are assumed to be
smooth enough (for all formulas to make sense) and linearly independent on $\T$.

\subsection{Prediction for  Markovian error processes}
\label{sec:mark}

\subsubsection{General Markovian process}

Consider the prediction of the random process \eqref{1.1} with $\T=[A,B]$ and the Markovian kernel $K(t,s)=u(t)v(s)$ for $t\le s$, where $u(\cdot)$ and $v(\cdot)$
are {twice differentiable} positive functions such that $q(t)=u(t)/v(t)$ is monotonically increasing.
As shown in \citep[Sect. 2.6]{DPZ2018}, a solution of the equation
$
 \int_A^B  K(t,s)\zeta(dt)= f(s)
$
holding for  all $s\in \T$
is the signed vector-measure $\zeta(dt)=z_{A}\delta_A(dt)+z_{B}\delta_B(dt)+z(t)dt$ with
\bea
 z_{A}&=&\frac{  1 }{ v^2(A) q^\prime (A)}  \Big[  \frac{f(A)u^\prime(A)}{u(A)} - f^\prime(A)  \Big]\, ,\;~~\\
 z(t)&=&-\frac{1}{ v( t )} \Big[ \frac{h^{\prime }(t)  }{q^{\prime } (t)}  \Big]^\prime\,,~~
 z_{B}=\frac{   h^\prime(B)}{ v(B) q^\prime (B)},
\eea
where $\psi'$ denotes  a derivative of a function $\psi$, the vector-function $h(\cdot)$ is defined by $h(t)=f(t)/v(t)$.

Then we obtain
$\zeta_{t_0}(dt)=z_{0A}\delta_A(dt)+z_{0B}\delta_B(dt)$ with
$$
 z_{0A}=\frac{  1-u'(A) }{ v^2(A) q^\prime (A)}v(t_0), ~~z_{0B}=\frac{1}{v(B)}v(t_0),
$$
$$
 C=\frac{1}{v^2(A)q(A)}f(A)f^\top (A)+\int_{A}^B \frac{[f(t)/v(t)]' [f(t)/v(t)]'^\top }{q'(t)}dt
$$
and
$$
 \tilde c=C^{-1}c=C^{-1}\left(f(t_0)-\left[\frac{  1-u'(A) }{ v^2(A) q^\prime (A)}v(t_0)f(A)+\frac{1}{v(B)}v(t_0)f(B)\right]\right).
$$
The BLUP measure is given by
$$
 Q_*(dt)=\zeta_{t_0}(dt)+c^\top G(dt)=\zeta_{t_0}(dt)+\tilde c^\top \zeta(dt)
$$
and the MSE of the BLUP is
$$
\mathrm{MSE}(\hat{y}_{Q_*}(t_0))=u(t_0)v(t_0)+\tilde c^\top   f(t_0)- \int_{A}^B K(t,t_0) Q_*(dt)\, .
$$

\subsubsection{Prediction when the error process is Brownian motion}

The covariance kernel $K(t,s)=\min(t,s)$ of Brownian motion is a particular case of the Markovian kernel with $u(t)=t$ and $v(s)=1$, $t\le s$.
Further we present the BLUP for few choices of $f(t)$.

For the location-scale model with $f(t)=1$, we obtain $c=0$ and, therefore, the BLUP measure is given by $Q_{*}(dt)=\delta_B(dt)$.
The BLUP  is $\hat y_{Q_*}(t_0)=y(B)$ and it has $\mathrm{MSE}(\hat y_{Q_*}(t_0))=t_0-B$.

For the  model with $f(t)=t$, we obtain $\tilde c=B^{-1}(t_0-B)$ and, thus, the BLUP measure is given by
$$
Q_{*}(dt)=\delta_B(dt)+\frac{t_0-B}{B}\delta_B(dt)=\frac{t_0}{B}\delta_B(dt).
$$
The BLUP  is $\hat y_{Q_*}(t_0)=\frac{t_0}{B} y(B)$ and it has
$\mathrm{MSE}(\hat y_{Q_*}(t_0))=\frac{t_0}{B} (t_0-B)$.

For the  model with $f(t)=t^2$, we obtain $\tilde c=(A^3+4/3(B^3-A^3))^{-1}(t_0^2-B^2)$ and,
thus, the BLUP measure is given by
$$
Q_{*}(dt)=\delta_B(dt)+\frac{t_0^2-B^2}{A^3+4/3(B^3-A^3)}\Big(2B\delta_B(dt)-A\delta_A(dt)-2dt\Big).
$$
The BLUP is
$$
\hat y_{Q_*}(t_0)= y(B)+\frac{t_0^2-B^2}{A^3+4/3(B^3-A^3)}\left(2By(B)-Ay(A)-2\int_A^B y(t)dt\right)
$$
and it has the mean squared error
$$
\mathrm{MSE}(\hat y_{Q_*}(t_0))=t_0+\tilde c \cdot t_0^2-\int_A^B t \cdot Q_{*}(dt).
$$

\subsubsection{Prediction for an OU error process}

The covariance kernel $K(t,s)=\exp(-\lambda|t-s|)$ of the OU error process is
also a particular case of the Markovian kernel with $u(t)=e^{\lambda t}$ and $v(s)=e^{-\lambda s}$, $t\le s$.

For the location-scale model $f(t)=1$, we obtain $\tilde c=(1+(B-A)\frac{\lambda}{2})^{-1}(1-e^{-\lambda|t_0-B|})$ and, therefore, the BLUP measure is given by
$$
Q_{*}(dt)=\tilde c/2\delta_A(dt)+(e^{-\lambda|t_0-B|}+\tilde c/2)\delta_B(dt)+\tilde c\lambda/2 dt.
$$
The BLUP is
$
\hat y_{Q_*}(t_0)=\tilde c/2 y(A)+(e^{-\lambda|t_0-B|}+\tilde c/2)y(B)+\tilde c\lambda/2 \int_A^B y(t)dt
$
and it has $\mathrm{MSE}(\hat y_{Q_*}(t_0))=1+\tilde c-\int_{A}^B e^{-\lambda|t-t_0|}Q_{*}(dt)$.

In Table~\ref{tab4} we give values of the square root of the MSE of the BLUP at the point $t_0=2$ for
the $N$-point equidistant design in the location scale model on the interval $[0,1]$
and the OU kernel with $\lambda=2$. From this table, we can see that one does not need many points to get almost optimal prediction: indeed, the MSE for designs with $N\ge 4$ is very close to the MSE for the continuous design. Similar results have been observed for other points $t_0$ and other Markovian kernels.

\begin{table}[!hhh]
\caption{The square root of the MSE of the BLUP at the point $t_0=2$ for
the $N$-point equidistant design in the location scale model on the interval $[0,1]$
and the OU kernel with $\lambda=2$; for the continuous design $\sqrt{\mathrm{MSE}}=1.164262$.}
\begin{center}
\begin{tabular}{|c|c|c|c|c|c|c|}
\hline
 $N$ & 2 & 4 & 8 & 16 &32 \\
\hline
$\sqrt{\mathrm{MSE}}$&1.18579 &1.167157 &1.164806 &1.164381& 1.16429 \\
\hline
\end{tabular}
\end{center}
\label{tab4}
\end{table}

\subsection{Prediction when the error process is integrated Brownian motion}
\label{sec:intB}

Consider the prediction of the random process \eqref{1.1} with $\T=[A,B]$, the 4 times differentiable vector of regression functions $f(t)$ and
the kernel of the integrated Brownian motion defined by
$$
K(t,s)={\min(t,s)^2}(3\max(t,s)-\min(t,s))/{6}.
$$
From \citep[Sect. 3.2]{DPZ2018} we have that
the signed matrix-measure $\zeta(dt)=(\zeta_0(dt),\zeta_1(dt))$ has
components $\zeta_0(dt)=z_A\delta_A(dt)+z_B\delta_B(dt)+z(t)dt$ and
$\zeta_1(dt)=z_{1,A}\delta_A(dt)+z_{1,B}\delta_B(dt)$, where
\bea
 z_A&=&f^{(3)}(A)-\frac{6}{A^2}f^{(1)}(A)+\frac{12}{A^3}f(A),\nonumber\\
 z_{1,A}&=&-f^{(2)}(A)+\frac{4}{A}f^{(1)}(A)-\frac{6}{A^2}f(A),\\
  z_B&=&-f^{(3)}(B),~~~ z_{1,B}=f^{(2)}(B),\;\;z(t)=f^{(4)}(t).\nonumber
\eea
Then we obtain
$\zeta_{t_0,0}(dt)=z_{t_0,A}\delta_A(dt)+z_{t_0,B}\delta_B(dt)+z_{t_0}(t)dt$ and $\zeta_{t_0,1}(dt)=z_{t_0,1,A}\delta_A(dt)+z_{t_0,1,B}\delta_B(dt)$ with (for $t_0>B$)
\bea
 z_{t_0,A}&=&K^{(3)}(A,t_0)-\frac{6}{A^2}K^{(1)}(A,t_0)+\frac{12}{A^3}K(A,t_0)=0,\\
 z_{t_0,1,A}&=&-K^{(2)}(A,t_0)+\frac{4}{A}K^{(1)}(A,t_0)-\frac{6}{A^2}K(A,t_0)=0,\\
  z_{t_0,B}&=&-K^{(3)}(B,t_0)=1,~z_{t_0,1,B}=K^{(2)}(B,t_0)=t_0-B,
\eea
and $z_{t_0}(t)=K^{(4)}(t,t_0)=0.$ This implies $\zeta_{t_0,0}(dt)=\delta_B(dt)$ and $\zeta_{t_0,1}(dt)=(t_0-B)\delta_B(dt)$.
Also we obtain
$
C=\frac{12}{A^3}f(A)f^\top (A)-\frac{6}{A^2}\Big(f'(A)f^\top (A)+f(A)f'^\top (A)\Big)+\frac{4}{A}f'(A)f'^\top (A)+\int_A^B f''(t)f''^\top (t)dt
$
and
$
\tilde c_0\!=\!C^{-1}c_0\!=\! C^{-1} \Big(f(t_0) -[f(B)+(t_0-B)f'(B)]\Big).
$


\vspace{2mm}
For the location-scale model with $f(t)=1$, we obtain $c=0$ and, therefore, the BLUP measure is given by
$\mathbf{Q}_{*}(dt)=(\delta_B(dt),(t_0-B)\delta_B(dt))^\top $.
The BLUP  is $\hat y_{\mathbf{Q}_*}(t_0)=y(B)+(t_0-B)y'(B)$ and it has
$\mathrm{MSE}(\hat y_{\mathbf{Q}_*}(t_0))=t_0^3/3-t_0B(t_0-B/2).$

\bigskip
\medskip

\noindent 	
{\bf Acknowledgments.}
This work has been supported in part by the Collaborative Research Center ``Statistical modelling of nonlinear
dynamic processes'' (SFB 823, Teilprojekt  C2) of the German Research Foundation (DFG).  The authors are grateful  to Martina Stein, who typed parts of this paper with considerable technical expertise. 

\bibliographystyle{elsarticle-harv}

\begin{thebibliography}{13}
\expandafter\ifx\csname natexlab\endcsname\relax\def\natexlab#1{#1}\fi
\expandafter\ifx\csname url\endcsname\relax
  \def\url#1{\texttt{#1}}\fi
\expandafter\ifx\csname urlprefix\endcsname\relax\def\urlprefix{URL }\fi

\bibitem[{Cressie(1993)}]{cressie1993statistics}
Cressie, N., 1993. Statistics for Spatial Data. John Wiley \& Sons.

\bibitem[{Dette et~al.(2019)Dette, Pepelyshev, and Zhigljavsky}]{DPZ2018}
Dette, H., Pepelyshev, A., Zhigljavsky, A., 2019. The {\sc blue} in
  continuous-time regression models with correlated errors. Annals of
  Statistics 47, 1928--1959.

\bibitem[{Fuentes(2006)}]{fuentes2006}
Fuentes, M., 2006. Testing for separability of spatial-temporal covariance
  functions. Journal of Statistical Planning and Inference 136~(2), 447--466.

\bibitem[{Gneiting et~al.(2007)Gneiting, Genton, and
  Guttorp}]{Gneitingetal2006}
Gneiting, T., Genton, M., Guttorp, P., 2007. Geostatistical space-time models,
  stationarity, separability and full symmetry. In: B.~Finkenstadt, L.~Held,
  V.~I. (Ed.), Statistical Methods for Spatio-temporal Systems. Chapman and
  Hall / CRC, Boca Raton, FL, pp. 151--176.

\bibitem[{Leatherman et~al.(2017)Leatherman, Dean, and
  Santner}]{leatherman2017designing}
Leatherman, E.~R., Dean, A.~M., Santner, T.~J., 2017. Designing combined
  physical and computer experiments to maximize prediction accuracy.
  Computational Statistics \& Data Analysis 113, 346--362.

\bibitem[{Morris et~al.(1993)Morris, Mitchell, and
  Ylvisaker}]{morris1993bayesian}
Morris, M.~D., Mitchell, T.~J., Ylvisaker, D., 1993. Bayesian design and
  analysis of computer experiments: use of derivatives in surface prediction.
  Technometrics 35~(3), 243--255.

\bibitem[{N{\"a}ther and {\v{S}}im{\'a}k(2003)}]{nather2003effective}
N{\"a}ther, W., {\v{S}}im{\'a}k, J., 2003. Effective observation of random
  processes using derivatives. Metrika 58~(1), 71--84.

\bibitem[{Parzen(1961)}]{parzen1961approach}
Parzen, E., 1961. An approach to time series analysis. The Annals of
  Mathematical Statistics 32~(4), 951--989.

\bibitem[{Rasmussen and Williams(2006)}]{rasmussen2006gaussian}
Rasmussen, C., Williams, C., 2006. Gaussian Processes for Machine Learning. MIT
  Press.

\bibitem[{Ripley(1991)}]{ripley1991statistical}
Ripley, B.~D., 1991. Statistical Inference for Spatial Processes. Cambridge
  University Press.

\bibitem[{Sacks et~al.(1989)Sacks, Welch, Mitchell, and Wynn}]{sacks1989design}
Sacks, J., Welch, W.~J., Mitchell, T.~J., Wynn, H.~P., 1989. Design and
  analysis of computer experiments. Statistical Science 4, 409--423.

\bibitem[{Santner et~al.(2003)Santner, Williams, and Notz}]{santner2003design}
Santner, T.~J., Williams, B.~J., Notz, W.~I., 2003. The Design and Analysis of
  Computer Experiments. Springer Series in Statistics. New York:
  Springer-Verlag.

\bibitem[{Stein(1999)}]{stein1999interpolation}
Stein, M.~L., 1999. Interpolation of Spatial Data: Some Theory for Kriging.
  Springer Science \& Business Media.

\end{thebibliography}

\end{document}